\documentclass[preprint,12pt]{elsarticle}

\usepackage{amssymb}
\usepackage{amsmath}
\usepackage{amsthm}
\usepackage[ruled,vlined]{algorithm2e}
\usepackage{graphicx}
\usepackage{textcomp}
\usepackage{array}
\usepackage{makecell}
\usepackage{booktabs}
\usepackage{longtable}
\usepackage[hidelinks]{hyperref}
\usepackage{dblfloatfix}
\usepackage{soulutf8}
\usepackage[utf8]{inputenc}
\usepackage{csquotes}
\usepackage[english]{babel}

\newtheorem{theorem}{Theorem}[section]

\newtheorem{lemma}[theorem]{Lemma}
\newtheorem{definition}{Definition}[section]
\newtheorem{proposition}{Proposition}[section]

\makeatletter

\let\c@algorithm\relax
\makeatother

\journal{European Journal of Operational Research (EJOR)}

\begin{document}

\begin{frontmatter}



\title{A Global Solution Algorithm for AC Optimal Power Flow through Linear Constrained Quadratic Programming}


\author[inst1]{Masoud Barati}

\affiliation[inst1]{organization={Electrical and Computer Engineering, University of Pittsburgh},
            addressline={Swanson School of Engineering}, 
            city={Pittsburgh},
            postcode={15206}, 
            state={PA},
            country={USA}}

\begin{abstract}
We formulate the Alternating Current Optimal Power Flow Problem (ACOPF) as a Linear Constrained Quadratic Program (LCQP) with a substantial number of negative eigenvalues ($r$) subjected to linear constraints, a problem acknowledged as NP-hard. We introduce two algorithms, Feasible Successive Linear Programming (FSLP) and Feasible Branch-and-Bound algorithm (FBB), aimed at a global optimal solution methodology. These algorithms employ a diverse range of simple but efficient optimization strategies such as bounded successive linear programming, convex relaxation, initialization, and branch-and-bound. Their objective is to discover a globally optimal solution to the inherent ACOPF, adhering to a predefined $\epsilon$-tolerance level. The complexity bound of the FSLP and FBB algorithms is $\mathcal{O}\left(N \prod_{i=1}^r\left\lceil\frac{\sqrt{r}\left(t_u^i-t_l^i\right)}{2 \sqrt{\epsilon}}\right\rceil\right)$, here, $N$ signifies the complexity involved in resolving the subproblems at each node within the framework of the FBB algorithm. The variables $t_l$ and $t_u$, elements of the set of real numbers raised to the power of $r$, represent the lower and upper boundaries of $t$ respectively. Additionally, $-|t|^2$ is identified as the only negative quadratic component present in the objective function of the elevated ACOPF problem. Second, we employ a combination of penalized semidefinite modeling, convex relaxation, and line search techniques to design a globally feasible branch-and-bound algorithm for the LCQP form of ACOPF, capable of locating an optimal solution to ACOPF within $\epsilon$-tolerance.
Our initial numerical results illustrate that the FSLP and FBB algorithms can remarkably identify a global optimal solution for large-scale ACOPF instances, even when $r$ is large. Furthermore, the suggested algorithm demonstrates enhanced performance over other cutting-edge methods in the majority of PG-lib test scenarios. 
\end{abstract}



\begin{keyword}
AC optimal power flow \sep brutal force algorithm \sep Branch and bound algorithm \sep global solution \sep lifted problem \sep linear constrained quadratic programming \sep nonconvex optimization \sep second order cone programming \sep successive linear program \sep semidefinite programming. \\

\hspace{0.8cm}

\small \textit{This work was supported by the National Science Foundation under Grant ECCS no. 1711921.}



\end{keyword}

\end{frontmatter}



\section{Introduction}
\label{sec:introduction}

The idea of AC Optimal Power Flow (ACOPF) problems was first introduced by Carpentier in the 1960s. The focus was to determine the optimal operational states for multiple generators within a transmission network, accommodating spatially and temporally varying demand \cite{1,2}. Rooted in the research surrounding power systems and their operation, OPF has matured into a well-founded discipline, utilized worldwide for day-to-day power grid management and regulation.


The ACOPF is an optimization problem widely encountered in different segments of power systems, strictly adhering to the AC power flow equations. This problem is usually structured to minimize the cost of generation while concurrently maintaining adherence to physical principles and engineering limitations. Conventional ACOPF problems are classified as Quadratically Constrained Quadratic Programs (QCQPs), which are complex mathematical models. According to the nonlinear characteristics of power flow equations, ACOPF problems display nonconvexity, rendering optimal power flow (OPF) problems intrinsically non-deterministic polynomial-time hard (NP-hard), as stated in the citation \cite{3}. Given their nonlinear and nonconvex characteristics, ACOPF problems have been a focal point of extensive research across multiple disciplines. Existing solutions for ACOPF problems are generally divided into two primary categories: those relying on mathematical algorithms and those employing ACOPF-specific algorithms.

From a mathematical perspective, ACOPF problems are characterized as quadratic programs (QPs) with a handful of negative eigenvalues, subjected to quadratic constraints. The extensive application scope of non-convex QPs, coupled with their computational complexity, has drawn the attention of researchers from various fields, resulting in a plethora of algorithms presented in the literature. Initial research efforts on non-convex QPs were predominantly concentrated on acquiring first-order or (weak) second-order Karush-Kuhn-Tucker (KKT) critical points with favorable objective values. For example, in reference \cite{4}, an uncomplicated alternating updating scheme was proposed, specifically designed for a QPs subclass, the bi-linear program. This method progressively converges to a KKT point. The authors expanded on this by detailing a local search process founded on linear optimization techniques to pinpoint a local optimum. In \cite{5}, a successive linear optimization approach was devised to ascertain a KKT point of QPs. In a similar vein, \cite{6} offered an interior point method to determine a KKT point of the related QPs. 

Over the past two decades, semidefinite optimization programs (SDPs) have rapidly evolved and have proven to be highly effective in addressing nonconvex QPs by providing tight relaxations and approximate solutions. Reference \cite{7} explores relaxations for nonconvex QCQPs rooted in semidefinite programming (SDP) and the reformulation-linearization technique (RLT). This study demonstrates the combined application of SDP and RLT constraints yields bounds significantly superior to those produced when either technique is used independently.

The concept of successive semidefinite relaxations (SDR) of a binary polytope and the confirmation of the process's finite convergence were explored in \cite{8}. 
Reference \cite{9} combined positive polynomial representations as sums of squares with moments theory to develop successive SDRs. This hierarchy was demonstrated to exhibit finite convergent behavior under certain conditions. The topic of successive SDRs pertaining to generic QPs was covered in \cite{10,11}.

However, as highlighted in \cite{12,13}, the number of variables and constraints dramatically increases with the rising levels of hierarchy in these procedures, which renders these methods computationally non-scalable. To remedy this, an alternative approach of formulating a branch-and-bound (B\&B) for QPs based on SDRs was proposed in \cite{14,15}. The adoption of the entirely co-positive program, as documented in \cite{16}, served to refine the B\&B methodology employed for LCQPs.
 Despite these B\&B-based global solvers being capable of solving many small-scale LCQPs, they can surpass time/memory constraints and generate only a suboptimal solution for problem sizes ranging between $50$ and $100$, as observed in \cite{16}.

More recently, \cite{17} proposed a global algorithm, ADMBB, applicable for quadratic programs subjected to linear and convex quadratic constraints and possessing a few negative eigenvalues. This algorithm can pinpoint a globally optimal solution for the fundamental QP within a pre-determined epsilon-tolerance. Unfortunately, this methodology is not computationally practical for large-scale applications, and the presented algorithm does not consistently provide bounded solution results.

\begin{figure}[t]
    \centering
    \includegraphics[width= 0.8\columnwidth]{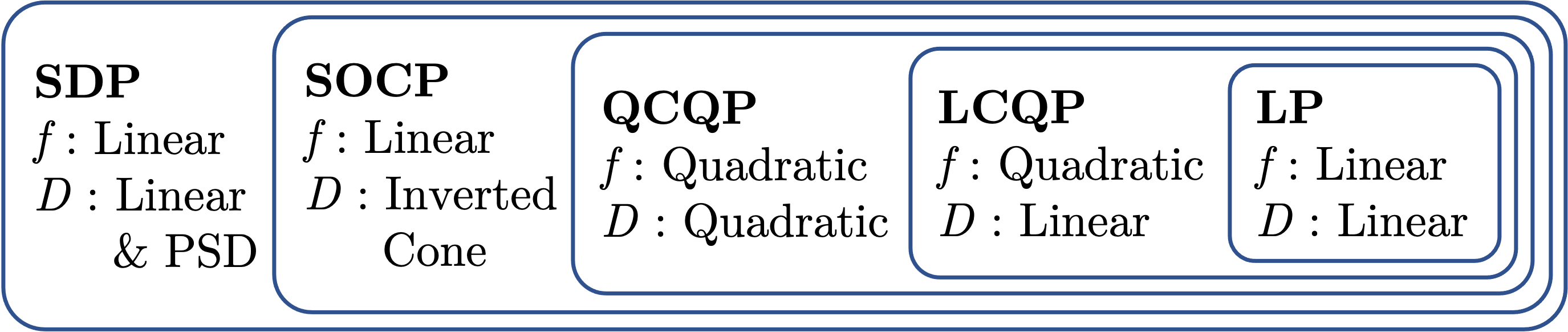}
    \caption{The atlas of ACOPF model in the literature including the specification of each model.}
    \label{fig: fig8}
\end{figure}

Due to the inherent nonconvexity of ACOPF problems, contemporary ACOPF-based algorithms are increasingly focusing on convexified methodologies. The methods mentioned can be broadly categorized into two distinct groups: relaxations and approximations, as stated in \cite{18}. The atlas of ACOPF models is depicted in Fig. \ref{fig: fig8}. As outlined in \eqref{eq: first}, the objective function $f$ can manifest as both linear and quadratic forms, while the feasible domain $\mathcal{D}$ can accommodate a range of constraint types, including linear, quadratic, inverted cone, and positive semidefinite linear constraints.
\begin{equation}{\label{eq: first}}
\min_{x \in \mathcal{D}} f(x)
\end{equation}

\vspace{-0mm}
Among the various algorithms, semi-definite programming (SDP) \cite{19,21} and second-order cone programming (SOCP) \cite{24,25} have been widely adopted and are well-established. However, it is crucial to mention that the conventional SOCP relaxation is not able to accommodate loop-based network constraints, specifically those where the total summation of voltage angle differences inside a loop equals zero. A contemporary investigation, referenced in \cite{28}, introduced a triad of enhancements to the second-order cone programming (SOCP) refinement applied to ACOPF. It is noteworthy, as discussed in \cite{26,28}, that the solution space permitted by the semidefinite programming (SDP) refinement is encompassed by that of the SOCP refinement.

Leveraging these techniques, a number of improved algorithms, such as convex iteration, have been proposed, showcasing superior efficiency \cite{20}. In \cite{21}, SDP was first used to address OPF problems, where OPF problems were reformulated into a SDP strcuture, and an interior point method (IPM) algorithm within SDP was developed. \cite{19} proposed a solution to the nonconvex ACOPF through the formulation of a convex semidefinite program based on the Lagrangian dual problem. The condition to ensure zero duality gap, derived from this work, was also satisfied for IEEE power test cases with up to 300 buses.

The research conducted in \cite{20} leveraged sparsity methodologies to improve the SOCP relaxation of ACOPF. We accomplished this by ensuring that the sub-matrices, associated with the largest cliques and smallest chord-free cycles in the cycle foundation, maintained a positive semidefinite state. Moreover, virtual lines were added to minimal chordless cycles, thereby transforming each into 3-node cycles. Additionally, a convex iteration was incorporated into their proposed solver, resulting in a convex relaxation that was tightly gapped.

The \cite{22} proposed a general solution framework for deriving near-optimal solutions for the SDP relaxation of QCQPs. They also introduced a rank penalty based on the nuclear norm of the Hermitian matrix to the SDP relaxation problem's objective, enforcing a low-rank solution. Reference \cite{29} proposed a penalized perturbed SDP relation with the aim to augment the real part of certain entries of the Hermitian matrix, thereby guaranteeing a rank-1 solution under the optimality condition. 
The larger the value of the regularization coefficient, the more substantial the penalty imposed on the rank of the Hermitian matrix.

The papers \cite{30,31} further investigated the presence of a rank-1 SDP solution in relation to obtaining a global solution for SCOPF with line contingencies. The authors developed a graph-theoretic convex program aimed at pinpointing network lines that pose problems, and integrated the loss incurred over these lines into the objective function, thus serving as a regularization term. Should the relaxation not be exact, penalizing the SDP problem serves as an alternative and feasible strategy.
 
In reference \cite{28}, a novel method utilizing convex quadratically constrained quadratic programs has been introduced. This method incorporates a penalty strategy that is consistent with semidefinite programming, second-order cone programming, and parabolic optimization techniques, ensuring the derivation of viable solutions for optimal power flow problems, given certain conditions. This method mainly tackles two prevalent challenges: firstly, the penalization constants are typically predefined using heuristic approaches; and secondly, the addition of regularization terms to the objective function in the ACOPF problem can alter its original value, as these terms persist even at the optimal solution.

Drawing inspiration from \cite{23} and the penalized SDP relaxation algorithm, this paper introduces two innovative algorithms: Feasible Successive Linear Programming (FSLP) and Feasible Branch-and-Bound algorithm (FBB). The principal contributions of this paper are summarized as follows:
\begin{enumerate}
    \item Elegant penalty terms, representing the second-order determinants of the symmetric semi-definite matrix, were integrated into the original objective function with the aim of achieving a rank of 1 for the matrix.
    \item The penalized SDP relaxation formulated in this study is standard quadratic programming with linear constraints (LCQP), making it simpler to solve compared to other formats.
    \item Taking into account the number of negative eigenvalues ($r$), two proposed algorithms were introduced: a bounded FSLP algorithm and a novel branch-and-bound (FBB) algorithm. These algorithms are designed to expedite the process of finding the optimal solution, surpassing the performance of current solution methodologies.
    \item When global optimal solutions are obtained, the penalty terms in the objective function are equal to zero. This indicates that these terms do not influence the objective values when optimal rank-1 solutions are acquired. This is a distinctive feature compared to the aforementioned methods, which used regularizers as penalty terms without exhibiting the vanishing property in the rank-1 solution result.
\end{enumerate}

The organization of this paper is outlined as follows: Section \uppercase\expandafter{\romannumeral2} presents the fundamental setup of the ACOPF formulation. Section \uppercase\expandafter{\romannumeral3} examines the penalized relaxation of SDP within the ACOPF problem. The different formats of ACOPF and the proposed bounded feasible successive linear programming algorithm are detailed in Section \uppercase\expandafter{\romannumeral4}. Section \uppercase\expandafter{\romannumeral5} delves into the proposed novel branch-and-bound algorithm, accompanied by supporting theorems and the corresponding proofs. The computational findings are detailed in Section \uppercase\expandafter{\romannumeral6}, while the conclusions are consolidated in Section \uppercase\expandafter{\romannumeral7}.

\section{AC Optimal Power Flow}
We delineate the construction of the ACOPF problem being examined within the context of this research work. 

\subsection{Mathematical ACOPF Formulation}
In the present study, we consider one of the most prevalent configurations of ACOPF problems in rectangular format, primarily centered on minimizing the generation expense while simultaneously adhering to physical laws and engineering restrictions. 
This section revisits the fundamental principles of AC power grid model in the form of optimization programming and integrates these principles to shed light on the pivotal ACOPF problem. It also establishes the notation that is consistently used throughout the paper.

Electric power grids incorporate different elements, including generating units, demands, nodes or buses, and transmission lines. In a broad analysis, the collections $\mathcal{N}$ of nodes and $\mathcal{L}$ of edges constitute a network graph denoted by $(\mathcal{N}, \mathcal{L})$, where the nodes correlate with buses and the edges with the connecting lines. It is essential to recognize $\mathcal{L}$ as a set of bidirectional edges. However, for each connecting line $(i, j) \in \mathcal{L}$, there is a designated directionality, with a ``source'' (from) side $(i, j)$ and a ``destination'' (to) side $(j, i)$.

\begin{figure}[t]
    \centering
    \includegraphics[width= 0.4\columnwidth]{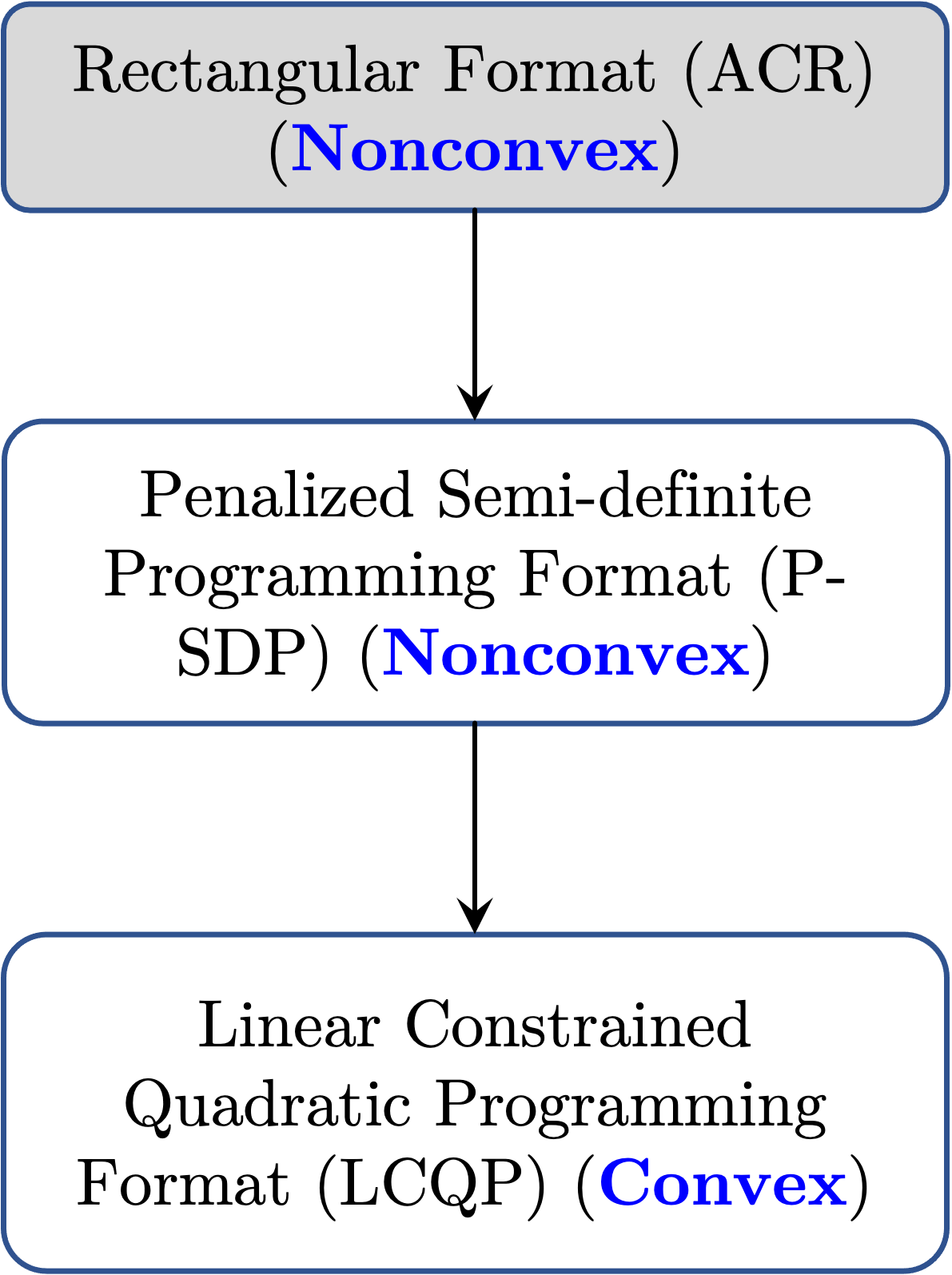}
    \caption{The trend of ACOPF model transformation in the presented solution methodology.}
    \label{fig: fig5}
\end{figure}
\vspace{-0mm}

thus enabling the tracking of line losses as power flows transition \emph{from-end} to \emph{to-end}. Usually, a slack bus $r \in \mathcal{N}$ is denoted for easy comparison of solutions and elimination of symmetric solutions.
AC power flow equations define the power flow within the network. These equations outline the relationship between complex quantities of current $I$, voltage $V$, power $S$, and admittance $y$, which are based on the physical principles of Kirchhoff's current law, Ohm's law, and the computation of the AC power grid. The integration of the aforementioned three principles results in the well-known network power flow equations,
\begin{align}
 S_i & = S_{G_i}-S_{D_i}=\sum_{(i,j) \in \mathcal{L}(i)} S_{ij} \quad \forall i \in \mathcal{N} \label{eq:s1}\\
 S_{ij} & =V_iI^*_{ij}=V_iy^*_{ij}(V^*_i-V^*_j) \nonumber \\
 & =y_{ij}^* V_i V_i^*-y_{i j}^* V_i V_j^* \quad \forall (i, j), \forall (j, i) \in \mathcal{L} \label{eq:s2}
\end{align}

Note that for bus $i \in \mathcal{N}$, the sum over $(i, j) \in \mathcal{L}$ integrates the power flow directed in the ``from'' orientation, while the sum over $(j, i) \in \mathcal{L}$ combines the flow directed in the ``to'' direction. Equations \eqref{eq:s1}-\eqref{eq:s2} constitute the basis for optimization of the power grid in operation and planning applications. These two equations are Nonconvex and nonlinear. However, each specific application enhances the functionality of these two equations with some additional supporting constraints. 
We employ a representative ACOPF problem formulation. This involves the conservation of active and reactive power flow at each node $i \in \mathcal{N}$, expressed as a function of the real and reactive power flow of each connecting line $(i,j) \in \mathcal{L}$ to every node $i \in \mathcal{N}$. The set of all lines connected to node $i$, is shown by $(i,j) \in \mathcal{L}(i)$. It also accounts for the maximum and minimum bounds of the real and reactive power flow of each connecting line $(i,j)$. Further considerations include upper and lower voltage magnitude limits at each node $i \in \mathcal{N}$, and the reference bus $r \in \mathcal{N}$ phase angle is set as zero. 
\begin{equation}{\label{eq: 3}}
\begin{split}
\text{min} & ~~ \sum_{i \in \mathcal{G}}~f_i(P_{Gi})   \\
s.t. &~ P_{ij}= (g_{i0}+g_{ij})(V_{di}^2+V_{qi}^2)-V_{di}(g_{ij}V_{dj}-b_{ij}V_{qj})  \\
  & \quad \ \ \ \ \  +V_{qi}(b_{ij}V_{dj}+g_{ij}V_{qj}), \qquad \
  \forall (i,j) \in \mathcal{L}   \\
  & Q_{ij}= -(b_{i0}+b_{ij})(V_{di}^2+V_{qi}^2)+V_{di}(b_{ij}V_{dj}+g_{ij}V_{qj}) \\
  & \quad \ \ \ \ \ \ \ \  -V_{qi}(g_{ij}V_{dj}-b_{ij}V_{qj}), \qquad \ \forall (i,j) \in \mathcal{L}   \\
  & ~   P_{G_i}-P_{D_i}=\sum_{(i,j) \in \mathcal{L}(i)} P_{ij}, \qquad \ \ \ \ \ \forall i \in \mathcal{N} \\
  & ~   Q_{G_i}-Q_{D_i}=\sum_{(i,j) \in \mathcal{L}(i)} Q_{ij}, \qquad \ \ \ \ \ \forall i \in \mathcal{N} \\
  & ~ -P_{ij}^{max} \leq P_{ij} \leq P_{ij}^{max}, \qquad \ \ \ \ \ \ \ \forall (i,j) \in \mathcal{L}   \\
  & ~ P_i^{min} \leq P_{G_i} \leq P_i^{max}, \qquad \ \ \ \ \ \ \ \ \ \ \ \ \forall i \in \mathcal{N}   \\
  & ~ Q_i^{min} \leq Q_{G_i} \leq Q_i^{max}, \qquad \ \ \ \ \ \ \ \ \ \ \ \ \forall i \in \mathcal{N}    \\
  & ~ (V_i^{min})^2 \leq V_{di}^2+V_{qi}^2 \leq (V_i^{max})^2, \quad \forall i \in \mathcal{N}  \\
  & ~ V_{qk} = 0. \quad   k \in \text{ref}
\end{split}
\end{equation}

The $\mathcal{G}$ represents the generator buses in \eqref{eq: 3}, while $\mathcal{N}=\{1,2,..., |\mathcal{N}|\}$ represents all buses within the system. The real and reactive power genartion units are denoted by $P_{G_i}$ and $Q_{G_i}$. The load is considered as a constant load with real and reactive power consumption $P_{D_i}$, and $Q_{D_i}$. The real and imaginary parts of the complex voltages $V_i$ are characterized as $V_{di}$ and $V_{qi}$. Furthermore, $g_{ij}$ and $b_{ij}$ represent the real and imaginary parts of the series admittance element of the line $(i,j) \in \mathcal{L}$ as functions of the series resistance $r_{ij}$ and reactance $x_{ij}$, where series conductance, $g_{ij}=\frac{r_{ij}}{r_{ij}^2+x_{ij}^2}$ and series susceptance, $b_{ij}=-\frac{x_{ij}}{r_{ij}^2+x_{ij}^2}$. The $g_{i0}=\frac{r_{i0}}{r_{i0}^2+x_{i0}^2}$ and $b_{i0}=-\frac{x_{i0}}{r_{i0}^2+x_{i0}^2}$ are the shunt elements connected to the bus $i \in \mathcal{N}$.
The quadratic monomials of the voltage terms featured in $P_{ij}$ and $Q_{ij}$ equations within \eqref{eq: 3} render the ACOPF problem nonconvex and categorize them as NP-hard. The \eqref{eq: 3} implies the ACOPF as a Quadratic Constrained Quadratic Program (QCQP) format. The presence of quadratic voltage terms in the active and reactive power flow equations, as seen in \eqref{eq: 3}, results in the ACOPF problems being nonconvex and thus classifies them as NP-hard. In response to the complexity and nonlinearity of this ACOPF problem, the penalized semidefinite programming algorithm (SDP) outlined in \cite{19} has been used, the details of which will be elaborated in the following section.

\begin{figure*}[t]
    \centering
    \includegraphics[width=1.0\linewidth]{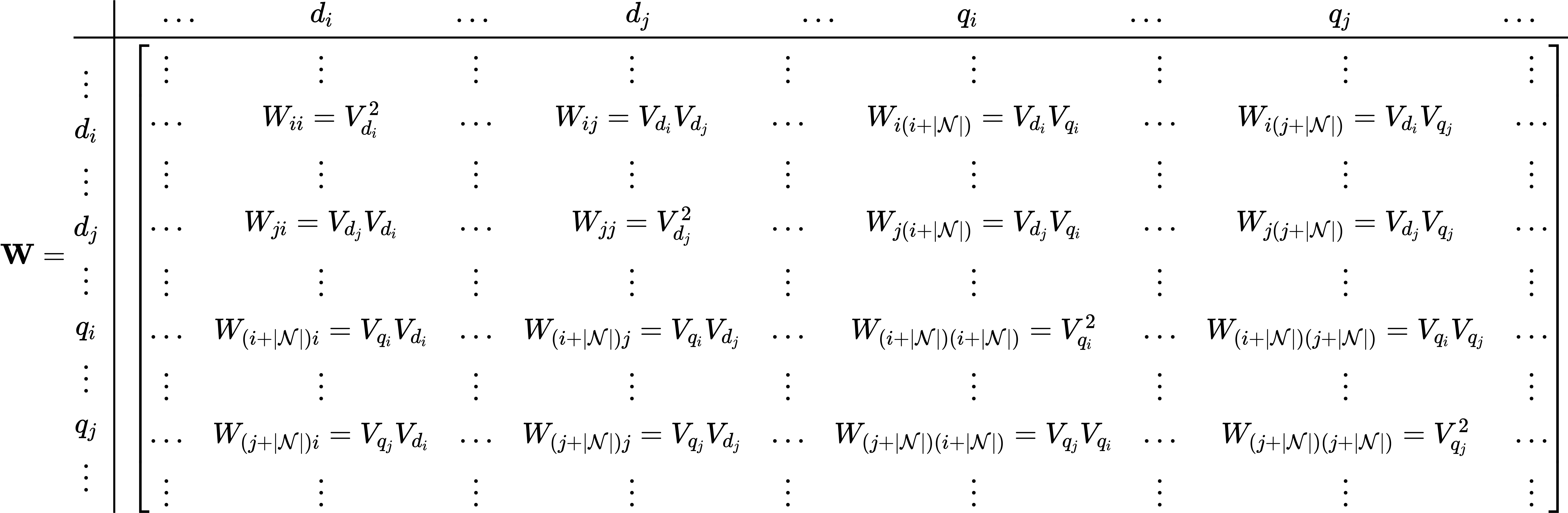}
    \caption{The shape of matrix \(\textbf{W}\) for two buses \(i\) and \(j\) includes all relevant elements in \(\textbf{W}\).}
    \label{fig:enter-label}
\end{figure*}


\section{Semidefine Programming (SDP) Relaxation}

This section initially involves the relaxation of the general ACOPF problem \eqref{eq: 3} utilizing semidefinite programming (SDP). To start with, we first define the notation. The vector $\mathbf{V}$ is the voltage vector of an $|\mathcal{N}|$-bus power system. 
\begin{equation}\label{eq: voltage}
\begin{aligned}
\mathbf{V}=[V_{d_1}, V_{d_2}, \dots , V_{d_{|\mathcal{N}|}} \mid V_{q_1}, V_{q_2}, \dots, V_{q_{|\mathcal{N}|}}]^{\top}
\end{aligned}
\end{equation}
Consider the $2|\mathcal{N}| \times 2|\mathcal{N}|$ Hermitian matrix defined as $\mathbf{W}$. 
Suppose $\mathbf{W}=\mathbf{V}\mathbf{V}^{\top}$ is defined as a $2|\mathcal{N}| \times 2|\mathcal{N}|$ symmetric matrix. In this case, equation \eqref{eq: 3} can be restructured utilizing the matrix $\mathbf{W}$ for two representative buses, $i$ and $j$, which correlate with respective constant voltage squares and bivariate voltage products. When considering a convex cost function, the semidefinite programming relaxation version can be derived from \eqref{eq: 3}, which is attained by dismissing the nonconvex rank-1 constraint in \eqref{eq: 4}.
The constraints characterizing $\mathbf{W}$ include:
\begin{equation}\label{eq: 4}
\begin{aligned}
\mathbf{W} & \succeq 0 \\
\operatorname{rank}(\mathbf{W}) & =1
\end{aligned}
\end{equation}
It should be emphasized that within the newly delineated W-space, the voltage constraints from \eqref{eq: 3} need to be presented in a different format.

\begin{equation}{\label{eq: 5}}
\begin{split}
\text{min} & ~ \sum_{i \in \mathcal{G}}~f_i(P_{Gi})  \\
s.t. &~ P_{ij}= (g_{i0}+g_{ij})(W_{ii}+W_{(i+|\mathcal{N}|)(i+|\mathcal{N}|)}) \\
  & \qquad  -g_{ij}W_{ij}+b_{ij}W_{i(j+|\mathcal{N}|)}-b_{ij}W_{(i+|\mathcal{N}|)j} \\
  & \qquad  -g_{ij}W_{(i+|\mathcal{N}|)(j+|\mathcal{N}|)}, 
  \qquad \forall (i,j) \in \mathcal{L}   \\
  & ~  Q_{ij}= -(b_{i0}+b_{ij})(W_{ii}+W_{(i+|\mathcal{N}|)(i+|\mathcal{N}|)})  \\
  & \qquad   +b_{ij}W_{ij}+g_{ij}W_{i(j+|\mathcal{N}|)}-g_{ij}W_{(i+|\mathcal{N}|)j}  \\
  & \qquad  +b_{ij}W_{(i+|\mathcal{N}|)(j+|\mathcal{N}|)}, 
  \qquad \forall (i,j) \in \mathcal{L}  
  \end{split}
\end{equation}
\begin{equation}
\begin{split}
  & ~  P_{G_i}-P_{D_i}=\sum_{(i,j) \in \mathcal{L}(i)} P_{ij}, \qquad \ \ \ \ \ \forall i \in \mathcal{N} \\
  & ~   Q_{G_i}-Q_{D_i}=\sum_{(i,j) \in \mathcal{L}(i)} Q_{ij}, \qquad \ \ \ \ \ \forall i \in \mathcal{N}
  \end{split}
\end{equation}
  \begin{equation}
\begin{split}
  & ~ -P_{ij}^{max} \leq P_{ij} \leq P_{ij}^{max}, \qquad \forall (i,j) \in \mathcal{L}   \\
  & ~ P_i^{min} \leq P_{G_i} \leq P_i^{max}, \qquad \forall i \in \mathcal{N}   \\
  & ~ Q_i^{min} \leq Q_{G_i} \leq Q_i^{max}, \qquad \forall i \in \mathcal{N}    
  \end{split}
\end{equation}
\begin{equation}
\begin{split}
  & ~ (V_i^{min})^2 \leq W_{ii}+W_{(i+|\mathcal{N}|)(i+|\mathcal{N}|)} \leq (V_i^{max})^2 \\
  & \qquad \qquad \qquad \qquad \qquad \ \ \ \ \ \ \ ~~~~~~~~~~\forall i \in \mathcal{N} \\
  & ~ W_{k(k+n)}=W_{(k+n)(k+n)}=0, \quad \exists~k \in \text{ref},  \\
  & ~ \operatorname{rank}(\mathbf{W})=1, \\
  & ~ \mathbf{W} \succeq 0.
\end{split}
\end{equation}

Presuming the generator cost functions are linear, i.e., $f_i(P_{G_i})=c_iP_{G_i}$, the objective function in \eqref{eq: 5} becomes a linear function of the $W$'s. 
In the redefined W-space, it is observed that the objective function assumes a linear form and is subjected to linear constraints. This SDP relaxation format yields a lower-bound solution that deviates significantly from an exact optimal solution. For an approximation that closely aligns with a near-global optimum, it is advisable to integrate a convexified variant of the Rank-1 constraint into the optimization formulation as expressed in equation \eqref{eq: 5}. In the forthcoming sections, we demonstrate that the Rank-1 constraint transforms into quadratic monomials in W-space within the objective function, which carries some negative eigenvalues— a characteristic indicative of a typical quadratic programming scenario. Consequently, we aim to resolve a quadratic program that presents a limited number of negative eigenvalues subject to linear constraints (LCQP). Within this framework, the symbol \(\operatorname{tr}\) denotes the trace operation on a matrix, defined as the aggregate of its diagonal entries. On the other hand, \(\succeq\) indicates that the matrix \(\mathbf{W}\) is positive semidefinite. $W_{ii}$ represents the $i$-th diagonal term in the $\mathbf{W}$ matrix, while the other forms of $W_{(\cdot)(\cdot)}$ in the Hermitian matrix are illustrated in the extensive matrix $\mathbf{W}$.
\begin{figure}[t]
    \centering
    \includegraphics[width= 0.7\columnwidth]{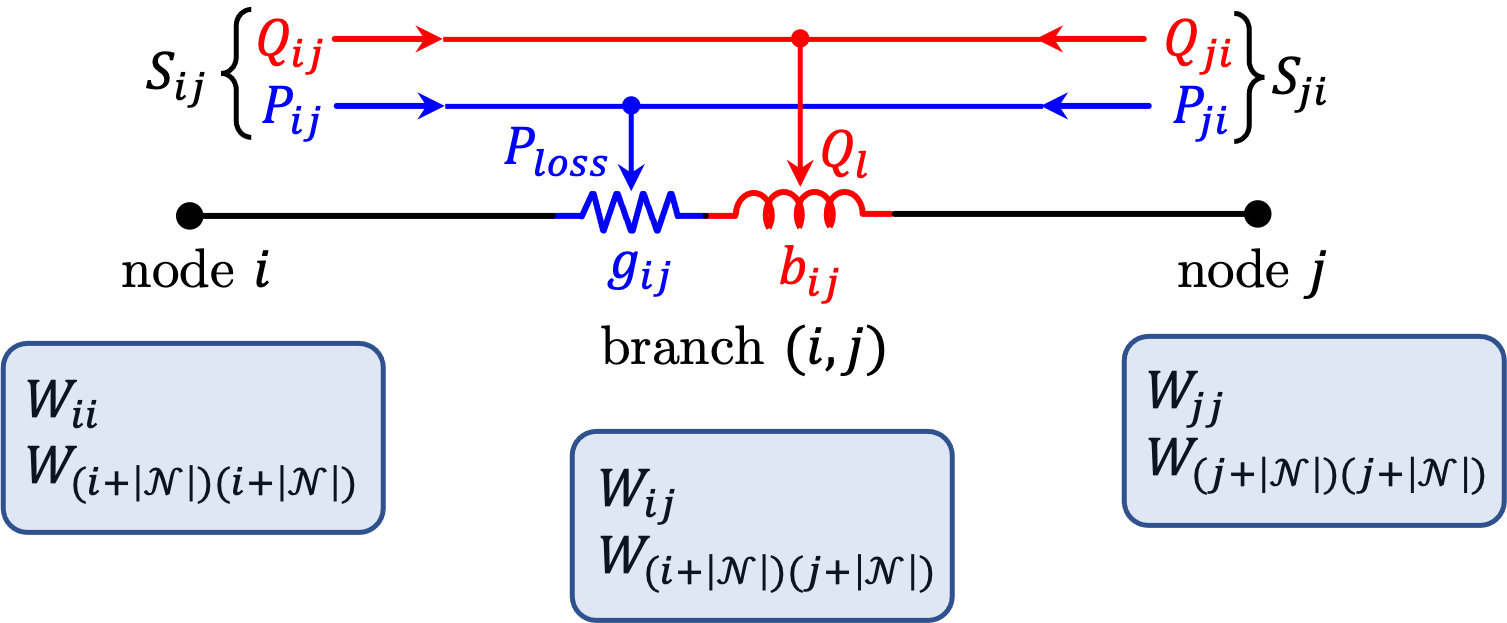}
    \caption{The corresponding elements of the Hermitian matrix related to nodes $i$ and $j$ of branch $(i,j)$ in \textbf{W}-space; the interplay between active and reactive power flow at the sending and receiving ends, taking into account ohmic losses and line reactive power consumption.}
    \label{fig: fig12}
\end{figure}
\vspace{-0mm}

Within this setting, solutions adhering to the $\operatorname{Rank}-1$ constraint are recognized as the globally optimal solutions to the ACOPF challenges.

\begin{theorem}\label{theorem: 1}
The $\operatorname{Rank}-1$ solutions for the matrix the second-order determinants of $W$ can reclaim $W$ due to the following statement. The $\operatorname{Rank}-1$ for the matrix $W$ implies that all $2 \times 2$ minor determinants of $W$ also possess a rank of $1$. \
\end{theorem}

\begin{proof}

Let's consider $\omega=\left(\begin{array}{cc}W_{ii} & W_{ij} \\ W_{ji} & W_{jj}\end{array}\right)$, which is a symmetric $2 \times 2$ matrix. Its leading principal minors are defined as $D_1=W_{ii}$ and $D_2=W_{ii} W_{jj}-W_{ij}^2$. To identify all the principal minors, we look at $\Delta_1=W_{ii}$ and $\Delta_1=W_{jj}$ (both of first order) as well as $\Delta_2=W_{ii} W_{jj}-W_{ij}^2$ (of second order).

We aim to understand the implications of leading principal minors being positive for the $2 \times 2$ matrices. Let's denote $\omega=\left(\begin{array}{ll}W_{ii} & W_{ij} \\ W_{ij} & W_{jj}\end{array}\right)$ as a symmetric $2 \times 2$ matrix. We need to demonstrate that if $D_1=W_{ii}>0$ and $D_2=W_{ii} W_{jj}-W_{ij}^2>0$, then $\omega$ is positively definite.

Given that $D_1=W_{ii}>0$ and $D_2=W_{ii} W_{jj}-W_{ij}^2>0$, we can also conclude that $W_{jj}>0$, since $W_{ii} W_{jj}>W_{ij}^2 \geq 0$ The characteristic equation of $\omega$ is
$$
\lambda^2-(W_{ii}+W_{jj}) \lambda+\left(W_{ii} W_{jj}-W_{ij}^2\right)=0
$$
and it gives us two solutions (as $A$ is symmetric), given by
$$
\lambda_{1,2}=\frac{W_{ii}+W_{jj}}{2} \pm \frac{\sqrt{(W_{ii}+W_{jj})^2-4\left(W_{ii} W_{jj}-W_{ij}^2\right)}}{2}
$$
As $(W_{ii}+W_{jj})>\sqrt{(W_{ii}+W_{jj})^2-4\left(W_{jj} W_{jj}-W_{ij}^2\right)}$, both solutions are positive, indicating that $\omega$ is positive definite.

At this juncture, we are prepared to verify this claim in a comprehensive way. Let's demonstrate this statement through mathematical induction. The base case where $n=1$ is quite straightforward. For the inductive step, let $\omega$ represent $2 \times 2$ principal submatrix of $W$ and define $d_m(t)=\operatorname{det}(\omega+tI)$. According to Jacobi's formula, $d_m^{\prime}(t)=\operatorname{tr}(\operatorname{adj}(\omega+t I))$. As the diagonal entries of adj$(\omega+t I)$ consist of some proper principal minors of $W+t I$, and considering that under the induction assumption, all proper principal submatrices of $W$ are positive semidefinite, $d_m^{\prime}(t)$ must be positive when $t>0$. This implies that $d_m$ strictly ascends to $[0, \infty)$. Therefore, for all $t>0$, $d_m(t)>d_m(0)=\operatorname{det}(\omega) \geq 0$.

This leads us to conclude that all principal minors of $W+t I$ are positive when $t>0$. Applying the original Sylvester's criterion for positive definite matrices, we find that $W+t I$ is positive definite. When taking $t \rightarrow 0^{+}$, it follows that $W$ is positive semidefinite.

Alternatively, the assertion can be easily proven without the need for mathematical induction. If we take $\omega$ to be any $2 \times 2$ principal submatrix of $W$ (this includes instances where $n=2$ and $\omega=W$) and $t>0$, then $\operatorname{det}(\omega+t I)=t^2+\sum_{k=1}^2 s_2(\omega) t^{2-k}$, where $s_2(\omega)$ represents the total summation of all the $2 \times 2$ principal minors of $\omega$. Each principal minor of $\omega$ also qualifies as a principal minor of $W, s_2(\omega)$ is non-negative. Thus, $\operatorname{det}(\omega+t I) \geq t^2>0$, implying that all principal minors of $W+t I$ are positive. Consequently, according to Sylvester's criterion for positive definite matrices, $W+t I$ is positive definite, which in turn makes $W=\lim _{t \rightarrow 0^{+}}(W+t I)$ positive semidefinite.
\end{proof}

\begin{theorem}\label{theorem: 2}
Furthermore, ohmic losses in the transmission line ensure that the cumulative sum of the two relevant minor determinants associated with the line $(i,j) \in \mathcal{L}$ is non-negative.
\end{theorem}
\begin{proof}
Ohmic losses in the transmission line connecting nodes $i$ and $j$ are quantified by the subsequent expression:
\begin{equation}
\begin{split}
&P_{ij}= (g_{i0}+g_{ij})(W_{ii}+W_{(i+|\mathcal{N}|)(i+|\mathcal{N}|)}) \\
  & \qquad  -g_{ij}W_{ij}+b_{ij}W_{i(j+|\mathcal{N}|)}-b_{ij}W_{(i+|\mathcal{N}|)j} \\
  & \qquad  -g_{ij}W_{(i+|\mathcal{N}|)(j+|\mathcal{N}|)}, 
\qquad \qquad \forall (i,j) \in \mathcal{L}   
\end{split}
\end{equation}
\begin{equation}
\begin{split}
&P_{ji}= (g_{j0}+g_{ji})(W_{jj}+W_{(j+|\mathcal{N}|)(j+|\mathcal{N}|)}) \\
& \qquad  -g_{ji}W_{ji}+b_{ji}W_{j(i+|\mathcal{N}|)}-b_{ij}W_{(j+|\mathcal{N}|)i} \\
& \qquad  -g_{ji}W_{(j+|\mathcal{N}|)(i+|\mathcal{N}|)}, 
\qquad  \qquad \forall (j,i) \in \mathcal{L}  
\end{split}
\end{equation}
\begin{equation}{\label{eq: ploss}}
\begin{split}
&P_{ij}+P_{ji}= g_{i0}W_{ii}+g_{j0}W_{jj}+g_{ij}(W_{ii}+W_{jj}-2W_{ij}) \\
& \qquad +g_{i0}W_{(i+|\mathcal{N}|)(i+|\mathcal{N}|)}+g_{j0}W_{(j+|\mathcal{N}|)(j+|\mathcal{N}|)} \\
& \qquad +g_{ij}(W_{(i+|\mathcal{N}|)(i+|\mathcal{N}|)}+W_{(j+|\mathcal{N}|)(j+|\mathcal{N}|)} \\
& \qquad \qquad \qquad -2W_{(i+|\mathcal{N}|)(j+|\mathcal{N}|)}),
\qquad \forall (i,j) \in \mathcal{L} 
\end{split}
\end{equation}
\begin{equation}
\begin{split}
&Q_{ij}= -(b_{i0}+b_{ij})(W_{ii}+W_{(i+|\mathcal{N}|)(i+|\mathcal{N}|)})  \\
& \qquad   +b_{ij}W_{ij}+g_{ij}W_{i(j+|\mathcal{N}|)}-g_{ij}W_{(i+|\mathcal{N}|)j}  \\
& \qquad  +b_{ij}W_{(i+|\mathcal{N}|)(j+|\mathcal{N}|)}, 
\qquad  \qquad \forall (i,j) \in \mathcal{L}  
\end{split}
\end{equation}
\begin{equation}
\begin{split}
& Q_{ji}= -(b_{j0}+b_{ji})(W_{jj}+W_{(j+|\mathcal{N}|)(j+|\mathcal{N}|)})  \\
& \qquad   +b_{ji}W_{ji}+g_{ji}W_{i(i+|\mathcal{N}|)}-g_{ji}W_{(j+|\mathcal{N}|)i}  \\
& \qquad  +b_{ji}W_{(j+|\mathcal{N}|)(i+|\mathcal{N}|)}, 
\qquad \forall (j,i) \in \mathcal{L}  
\end{split}
\end{equation}
\begin{equation}{\label{eq: qloss}}
\begin{split}
&Q_{ij}+Q_{ji}= -b_{i0}W_{ii}-b_{j0}W_{jj}-b_{ij}(W_{ii}+W_{jj}-2W_{ij}) \\
& \qquad  -b_{i0}W_{(i+|\mathcal{N}|)(i+|\mathcal{N}|)}-b_{j0}W_{(j+|\mathcal{N}|)(j+|\mathcal{N}|)} \\
& \qquad  -b_{ij}(W_{(i+|\mathcal{N}|)(i+|\mathcal{N}|)}+W_{(j+|\mathcal{N}|)(j+|\mathcal{N}|)} \\
& \qquad \qquad \qquad -2W_{(i+|\mathcal{N}|)(j+|\mathcal{N}|)}),
\qquad \forall (i,j) \in \mathcal{L} 
\end{split}
\end{equation}

The terms $W_{ii}$ and $W_{jj}$ maintain non-negative values. Additionally, according to Theorem \ref{theorem: 1}, we have $W_{ii}+W_{jj}-2W_{ij} \ge 0$ and $W_{(i+|\mathcal{N}|)(i+|\mathcal{N}|)}+W_{(j+|\mathcal{N}|)(j+|\mathcal{N}|)}-2W_{(i+|\mathcal{N}|)(j+|\mathcal{N}|)} \ge 0$. Consequently, we can infer that,   
\begin{equation}\label{eq: ploss2}
  P_{ij}+P_{ji} \ge 0,  
\end{equation}
Furthermore, since $b_{i0}$, $b_{j0}$, and $b_{ij}$ all have positive values, 
\begin{equation}\label{eq: qloss2}
  Q_{ij}+Q_{ji} \ge 0.  
\end{equation}
\end{proof}
Therefore, incorporating two fresh inequalities \eqref{eq: ploss2} and \eqref{eq: qloss2} into the SDP relaxation method can steer the solution towards achieving a precise result.

Finally, leveraging Theorem \ref{theorem: 1} and \ref{theorem: 2}, we introduce a penalized relaxation approach of SDP subject to linear constraints as a function of $W$s. In this methodology, the minor determinants are incorporated as new penalty terms in the objective function of \eqref{eq: 5}, along with the consideration of constraints \eqref{eq: ploss2} and \eqref{eq: qloss2}. This strategy ensures the tightest relaxation for SDP. Remarkably, this novel formulation can significantly steer the rank of the $W$ matrix towards $1$.
\subsection{Penalized SDP Relaxation}
The SDP relaxation with an imposed penalty for the equation denoted as \eqref{eq: 5} is expressed as follows,
\begin{equation} {\label{eq: 14}}
\begin{split}
  & \min \quad \sum_{k \in \mathcal{G}}f_k(P_{Gk})+\sum_{(i,j) \in \mathcal{L}} \lambda_{ij} (W_{ii}W_{jj}-W_{ij}^2) \\
  &~~~~~~ +\sum_{(i,j) \in \mathcal{L}} \lambda_{ij} (W_{(i+|\mathcal{N}|)(i+|\mathcal{N}|)}+W_{(j+|\mathcal{N}|)(j+|\mathcal{N}|)} \\
  & \qquad \qquad \qquad -2W_{(i+|\mathcal{N}|)(j+|\mathcal{N}|)}) \\
\end{split}
\end{equation}
\begin{equation}{\label{eq: 15}}
\begin{split}
s.t. &~ P_{ij}= (g_{i0}+g_{ij})(W_{ii}+W_{(i+|\mathcal{N}|)(i+|\mathcal{N}|)}) \\
  & \qquad  -g_{ij}W_{ij}+b_{ij}W_{i(j+|\mathcal{N}|)}-b_{ij}W_{(i+|\mathcal{N}|)j} \\
  & \qquad  -g_{ij}W_{(i+|\mathcal{N}|)(j+|\mathcal{N}|)}, 
  \qquad \forall (i,j) \in \mathcal{L}   \\
  & ~  Q_{ij}= -(b_{i0}+b_{ij})(W_{ii}+W_{(i+|\mathcal{N}|)(i+|\mathcal{N}|)})  \\
  & \qquad   +b_{ij}W_{ij}+g_{ij}W_{i(j+|\mathcal{N}|)}-g_{ij}W_{(i+|\mathcal{N}|)j}  \\
  & \qquad  +b_{ij}W_{(i+|\mathcal{N}|)(j+|\mathcal{N}|)}, 
  \qquad \forall (i,j) \in \mathcal{L}   \\
  & ~  P_{ij}+P_{ji} \ge 0,  \qquad \forall (i,j),(j,i) \in \mathcal{L} \\
  & ~  Q_{ij}+Q_{ji} \ge 0.  \qquad   \forall (i,j),(j,i) \in \mathcal{L} \\
  & ~ -P_{ij}^{max} \leq P_{ij} \leq P_{ij}^{max}, \qquad \forall (i,j) \in \mathcal{L} 
\end{split}
\end{equation}
\begin{equation}{\label{eq: 16}}
\begin{split}
  & ~  P_{G_i}=P_{D_i}+\sum_{(i,j) \in \mathcal{L}(i)} P_{ij}, \qquad \ \ \ \ \ \forall i \in \mathcal{N} \\
  & ~  Q_{G_i}=Q_{D_i}+\sum_{(i,j) \in \mathcal{L}(i)} Q_{ij}, \qquad \ \ \ \ \ \forall i \in \mathcal{N} \\
  & ~ P_i^{min} \leq P_{G_i} \leq P_i^{max}, \qquad \forall i \in \mathcal{N}   \\
  & ~ Q_i^{min} \leq Q_{G_i} \leq Q_i^{max}, \qquad \forall i \in \mathcal{N}    
\end{split}
\end{equation}
\begin{equation}{\label{eq: 17}}
\begin{split}
  & ~ (V_i^{min})^2 \leq W_{ii}+W_{(i+|\mathcal{N}|)(i+|\mathcal{N}|)} \leq (V_i^{max})^2 \\
  & \qquad \qquad \qquad \qquad \qquad \ \ \ \ \ \ \ ~~~~~~~~~~\forall i \in \mathcal{N} \\
  & ~ W_{k(k+n)}=W_{(k+n)(k+n)}=0. \quad \exists~k \in \text{ref}
\end{split}
\end{equation}
In equation \eqref{eq: 14}, the hyperparameter $\lambda_{ij} = \rho ((g_{i0}+g_{ij})^2+(b_{i0}+b_{ij})^2)$ signifies the relative importance of each transmission line in the objective function. Here, $\rho$ is a large value (for instance, 1000). The problem defined by equations \eqref{eq: 14} - \eqref{eq: 17} presents a nonconvex quadratic function with negative eigenvalues under linear constraints, constituting a canonical Linear Constrained Quadratic Programming problem (LCQP). The following section will detail the methodology to solve this LCQP problem using Feasible Successive Linear Programming ($\operatorname{FSLP}$) and a new branch and bound (BB) algorithm.

\subsection{ACOPF in the form of LCQP}
\begin{figure}[t]
    \centering
    \includegraphics[width= 1.0\columnwidth]{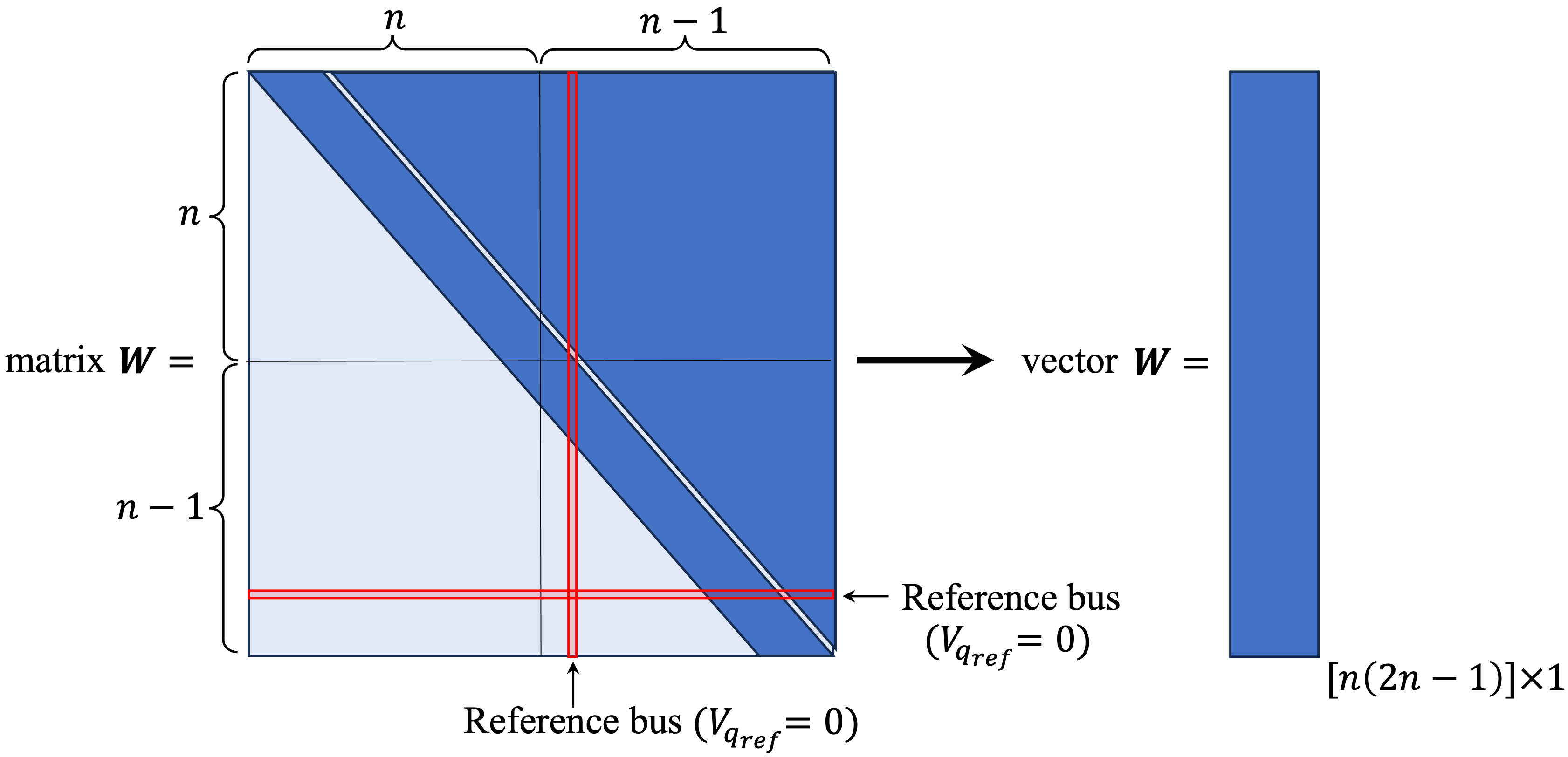}
    \caption{the shape of \textbf{W} matrix.}
    \label{fig: fig4}
\end{figure}

Firstly, we reformulate problem \eqref{eq: 14} - \eqref{eq: 17} into a standard quadratic form. This simplification will facilitate a more straightforward expression.
\begin{equation}\label{newpro}
\begin{split}
    & \min \quad f(W)=W^\top AW+b^\top W \\
    & \ s.t.\quad \ W \in \mathcal{F}, \ \mathcal{F}=\{W \in \mathbb{R}^n: LW \leq l\}
\end{split}
\end{equation}
In this case, $W$ denotes a vector in $\mathbb{R}^{n}$, $A$ signifies an $n\times n$ matrix in $\mathbb{R}$, $L$ represents an $m\times n$ matrix in $\mathbb{R}$, $b$ stands for a vector in $\mathbb{R}^n$, and $l$ is a vector in $\mathbb{R}^m$. Crucially, the $W$ vector in equation \eqref{newpro} is a reconfigured version of the $W$ matrix found in equations \eqref{eq: 14} - \eqref{eq: 17}, taking advantage of the symmetric attribute of the $W$ matrix. Figure \ref{fig: fig4} illustrates the new dimension of this vector. The matrix $W$ has active elements shown in blue that correspond to the components of the Hermitian matrix, equating to $(2|\mathcal{N}|-1)+\frac{(2|\mathcal{N}|-1)(2|\mathcal{N}|-1)-(2|\mathcal{N}|-1)}{2}$. Therefore, the size of the $W$ vector is $n \times 1$, where $n=(2\mathcal{N}-1)\mathcal{N}$. 
\vspace{-0mm}
Considering the presence of negative eigenvalues in matrix \(A\), the quadratic term in the objective function can be reformulated through the application of the Difference Convex (D.C.) program, as demonstrated in \cite{32},
\begin{equation}
\begin{split}
    W^\top AW & = W^\top A_{+}W-W^\top A_{-}W \\
          & = W^\top A_{+}W-\|CW\|^2,
\end{split}
\end{equation}
Matrices that are symmetric and exhibit non-negative definiteness \(A_+\) and \(A_-\) originate from the singular value decomposition (SVD) of \(A\). Moreover, \(C\) is an \(r \times n\) matrix, with \(r\) representing the count of negative eigenvalues in \(A\). Details regarding the matrix \(C\) can be found in the \emph{FBB Algorithm} section.
Hence, the problem \eqref{newpro} can be reframed as the following elevated optimization problem,
\begin{equation}\label{repro}
    \min_{(W,t)\in \mathcal{F}} b^\top W + W^\top A_{+}W - \|t\|^2,
\end{equation}
\begin{equation*}
    \mathcal{F} = \{(W,t) \mid LW \leq l, \ t_l \leq t \leq t_u, \ CW = t\}.
\end{equation*}

The notation $t_l,t_u \in \mathbb{R}^r$ signifies the respective lower and upper limits, which are obtained by solving the following linear optimization problems,
\begin{equation} \label{t}
    t_l^i=\min_{W\in \mathcal{F}} c_i^\top W,\ \ t_u^i=\max_{W \in \mathcal{F}} c_i^\top W, \ \ i=\{1,...,r\},
\end{equation}
Here, the term $c_i^\top$ denotes the transpose of the $i$-th row of the matrix $C$.
Similar to the approach used in \cite{33,34}, a linear surrogate objective function can be utilized to represent the negative quadratic component over a specific set $t$ bounded by the range $[t_l,t_u]$. Consequently, we derive the ensuing linear representation:
\begin{equation} \label{linear}
\begin{split}
    & \min \quad W^\top A_+W+b^\top W-t^\top CW, \\
    & \ s.t. \quad \ W\in \mathcal{F}.
\end{split}
\end{equation}
Expanding on the resolution of the issue delineated in \eqref{linear}, we have developed an algorithm for successive linear programming that is viable, denoted as $\operatorname{FSLP}$, to address the problem referenced in \eqref{repro}. This algorithm employs an iterative approach, sequentially updating $W$ and $t$ at each iteration to find the optimal solution. 

\section{Bi-convex format of the ACOPF and Feasible Successive Linear Programming Algorithm}

Subsequently, we examine a distinct Lagrangian approach for the problem denoted by \eqref{repro}, described in the following manner,
\begin{equation}
\mathcal{L}(W,t,\lambda)=W^\top A_+W + \lambda^\top (CW-t) + b^\top W - \|t\|^2,\label{L}
\end{equation}
In this scenario, the Lagrangian multipliers are utilized for the constraint $CW=t$. A thorough examination of the pertinent Karush–Kuhn-Tucker (KKT) conditions related to the previously mentioned Lagrangian function reveals $t$. It can be deduced that $\frac{\partial \mathcal{L}(W,t,\lambda)}{\partial t} = -2t - \lambda = 0$, thereby leading to the inference that $\lambda = -2t$. It implies that $\lambda$ is not an additional variable. 

It is important to note that the solution achieving optimality for the problem delineated in \eqref{newpro} is required to fulfill the criteria set by the KKT conditions associated with the respective Lagrangian function, as pointed out in \cite{35}. Consequently, the Lagrangian function can be reformulated in the following manner:
\begin{equation}
\mathcal{L}^*(W,t)=W^\top A_+W+b^\top W-2t^\top CW+\|t\|^2\label{re}
\end{equation}
For convenience, we call the above Lagrangian function as the optimal Lagrangian function. Because
\begin{equation*}
\begin{split}
    & W^\top A_+W+b^\top W-2t^\top CW+\|t\|^2 \\
     = \ & W^\top A_+W+b^\top W-\|CW\|^2+\|t-CW\|^2,
\end{split}
\end{equation*}
therefore, the problem \eqref{newpro} can be reformulated as a bi-convex optimization problem.
\begin{equation}\label{bi}
\begin{split}
& \min \quad W^\top A_{+}W+b^\top W-2t^\top CW+\|t\|^2 \\
& \ s.t. \quad\ W \in \mathcal{F}
\end{split}
\end{equation}
Although the method draws on traditional successive linear programming techniques, there is no guarantee in this particular problem that a bounded solution will be found. This means that there is a potential for $t_l^i$ to become $-\infty$ and $t_u^i$ to become $+\infty$. As such, our proposed $\operatorname{FSLP}$ offers an enhanced bounded successive linear programming strategy that avoids unbounded solutions. 
Using an elementary case of Young's inequality \cite{36}, the inequality with exponent 2 for two variables $t$ and $CW$ can lead to \eqref{Young}. 
\begin{equation}\label{Young}
\begin{split}
& t^\top CW \leq \frac{\|t\|^2}{2}+\frac{\|CW\|^2}{2}
\end{split}
\end{equation}
The inequality \eqref{Young} provides an upper bound limit on the bi-convex term and the following inequality provides a lower bound limit on the bi-convex term in the objective function. This inequality is based on the first-order taylor series on bi-variant term $t^\top (CW)$.
\begin{equation}\label{lower}
\begin{split}
& t^\top (CW) \ge t^{k^\top} C W^k + (t - t^k)^\top C W^k + t^{k^\top} C (W - W^k)
\end{split}
\end{equation}
Therefore, an enhanced bi-convex optimization problem can be written in \eqref{bi-new}.

\begin{equation}\label{bi-new}
\begin{split}
\min & \quad W^\top A_{+}W+b^\top W-2t^\top CW+\|t\|^2 \\
s.t. &  \quad t^\top CW \leq \frac{\|t\|^2}{2}+\frac{\|CW\|^2}{2} \\
& \quad t^\top CW \ge t^{k^\top} C W^k + (t - t^k)^\top C W^k \\
& \quad \quad \quad \quad \quad 
+ t^{k^\top} C (W - W^k) \\
& \quad W \in \mathcal{F}. 
\end{split}
\end{equation}
Notably, the strategy of transforming non-convex QP issues into bi-convex optimization is commonly found in literature \cite{37,38,39}. The unique aspect of the model \eqref{bi-new} is the bounded condition imposed on the variable $t$. The bi-convex optimization proposed in \cite{37,38,39} has no restrictions on $t$, leading to unbounded solution results. We can observe that the solution to the problem \eqref{bi-new} can enhance the value of the objective function more efficiently compared to its unrestricted counterpart. It is crucial to note that for a given $t$, the problems \eqref{bi-new} and \eqref{linear} are equivalent.

Given a choice of $t$ within the interval $[t_l,t_u]$, it becomes clear that the problem outlined in equation \eqref{bi-new} is feasible, bounded and properly formulated. Then, we can infer several important theorems.

\begin{theorem}
If $\Bar{W}$ is an optimal solution to the problem denoted by \eqref{bi-new}, under the condition that $t = C\Bar{W}$, then $\Bar{W}$ qualifies as a KKT point for \eqref{bi-new}.
\end{theorem}

Let \(\mathcal{S}(t)\) denote the set of optimal solutions for \eqref{bi-new}, and assume that \(W(t)\) is an element of \(\mathcal{S}(t)\). The following lemma outlines several essential characteristics of the mapping \(CW(t)\) in relation to \(t\).

\begin{lemma}
\label{pro32}
Given that \(W(t)\in \mathcal{S}(t)\), we observe the following:
\begin{enumerate}
    \item For any \(W \in \mathcal{F}\) satisfying the inequality
    \begin{equation*}
    |CW-t| - |CW(t)-t|\leq 0,
    \end{equation*}
    it follows that
    \begin{equation*}
    f(W) - f(W(t)) \geq 0,
    \end{equation*}
    
    \item If \(W^*\) represents the global optimal solution of \eqref{bi-new}, then
    \begin{equation}
    |CW(t)-t| - |CW^*-t|\leq 0,
    \end{equation}
    
    \item The function \(CW(t)\) is monotonic with respect to \(t\), as evidenced by
    \begin{equation} \label{s3}
    (t_1-t_2)^\top [CW(t_1)-CW(t_2)] \geq 0,
    \end{equation}
    where
    \begin{equation*}
    t_1 \neq t_2 \in [t_l,t_u],\ W(t_1) \in \mathcal{S}(t_1),\ W(t_2) \in \mathcal{S}(t_2).
    \end{equation*}
\end{enumerate}
\end{lemma}

\begin{proof}
The initial assertion is a direct consequence of the recognition that \(W(t)\) constitutes the optimal solution for problem \eqref{bi}, coupled with the realization that \(|CW-t| \leq |CW(t)-t|\). Now, let us turn our attention to the second proposition. Given that \(W^*\) represents the global optimum of the problem and considering \(W(t) \in \mathcal{F}\), we can confidently state that
\begin{equation} \label{p1}
f(W(t)) - f(W^*) \geq 0.
\end{equation}
Given that $W^*$ is a feasible solution for problem \eqref{bi}, and given $W(t)\in \mathcal{S}(t)$, it follows that
\begin{equation} \label{p2}
f(W(t))+\|CW(t)-t\|^2 - f(W^*) - \|CW^*-t\|^2 \leq 0.
\end{equation}
From equations \eqref{p1} and \eqref{p2}, it is evident that
\begin{equation}
\|CW(t)-t\| - \|CW^*-t\|\leq 0.
\end{equation}
This verifies second statement of the theorem.\\
\indent Next, for third statement, assume $t_1 \neq t_2 \in [t_l,t_u]$, $W(t_1) \in \mathcal{S}(t_1)$, and $W(t_2) \in \mathcal{S}(t_2)$. Using the optimal conditions of problem \eqref{bi}, we can deduce that
\begin{align*}
& f(W(t_1))+\|CW(t_1)-t_1\|^2 \leq f(W(t_2))+\|CW(t_2)-t_1\|^2, \\
& f(W(t_2))+\|CW(t_2)-t_2\|^2 \leq f(W(t_1))+\|CW(t_1)-t_2\|^2.
\end{align*}
Adding the above two inequalities, we derive
\begin{align*}
& t_2^\top CW(t_2)+t_1^\top CW(t_1) \geq t_2^\top CW(t_1)+t_1^\top CW(t_2) \\
& \Longleftrightarrow (t_1-t_2)^\top C(W(t_1)-W(t_2)) \geq 0,
\end{align*}
which yields relation \eqref{s3}. This concludes the proof.

Moreover, the lower bound constraint \eqref{lower} provides a comprehensive way to prove the statement. The linearized bi-variant function at $t_1$ can be expressed as follows.
\begin{equation*}
\begin{split}
& t^\top CW(t)\mid_{t=t_1} \ge t_1^\top CW(t_1) + (t-t_1)^\top CW(t_1) \\
& \qquad \qquad \qquad + t_1^\top C(W(t)-W(t_1)) \\
\end{split}
\end{equation*}
Hence, at $t=t_2$, we have, 
\begin{equation*}
\begin{split}
& t_2^\top CW(t_2) \ge t_1^\top CW(t_1) + (t_2-t_1)^\top CW(t_1) \\
& \qquad \qquad \qquad + t_1^\top C(W(t_2)-W(t_1)) \\
\end{split}
\end{equation*}
It implies that, 
\begin{equation*}
\begin{split}
& (t_2-t_1)^\top CW(t_2) \ge (t_2-t_1)^\top CW(t_1) \\
& (t_1-t_2)^\top CW(t_2) \leq (t_1-t_2)^\top CW(t_1) 
\end{split}
\end{equation*}
Finally, 
\begin{equation*}
\begin{split}
& (t_1-t_2)^\top C(W(t_1)-W(t_2)) \ge 0 
\end{split}
\end{equation*}
\end{proof}

The method of Feasible Successive Linear Programming (FSLP) utilized to tackle LCQP problems is encapsulated in the algorithm delineated below. 

\begin{algorithm}[H]
\caption{Feasible Successive Linear Programming $\operatorname{FSLP}$($W^0,\epsilon$)}
\label{a1}
\KwIn{Commencing with \(W^0\) as the starting value and employing \(\epsilon > 0\) as the threshold for termination}
\KwOut{The ultimate solution $(W^{k+1},t^{k+1})$ is returned when the algorithm halts}

\textbf{Step 0} \textit{Initialization:} Set the iteration counter $k$ to 0. Calculate the initial value of $t^k$ using the formula $t^k=CW^k$, where $C$ is a given matrix, and $W^k$ is the current estimate of the parameter values\;

\textbf{Step 1} \textit{Iterative Optimization:} Solve the given optimization problem (denoted by equation \eqref{bi-new}) with $t$ set to the current $t^k$ to obtain the optimal solution $W(t^k)$. Update the parameter estimates by setting $W^{k+1} = W(t^k)$. Calculate the new value of $t$ for the next iteration using the formula $t^{k+1}=CW^{k+1}$\;

\textbf{Step 2} \textit{Convergence Check:} Compute the Euclidean distance between $t^{k+1}$ and $t^k$. If this distance exceeds the square root of the stopping criterion $\sqrt{\epsilon}$, increase the number of iterations $k$ by 1 and return to the step of ``Iterative Optimization''. If this condition does not hold, assume that the solution has converged\;

\end{algorithm}

We need supplementary supporting properties to ensure that the solution of \eqref{bi-new} behaves appropriately throughout the $\operatorname{FSLP}$ iterative sequence. The following Lemma offers proof of the solution's monotonic behavior at each step in the iteration.  

\begin{lemma}\label{le31}
Consider the sequences \({W^k}\) and \({t^k}\), produced through the \(\operatorname{FSLP}\) algorithm, with \(W^*\) representing the global optimal solution of \eqref{newpro}. The following propositions are valid:

\eqref{newpro}. The following assertions hold true:
\begin{enumerate}
    \item The sequence ${t^k}$ exhibits a monotonic trend, as can be seen from
    \begin{equation*}
    (t^k - t^{k-1})^\top (t^{k+1} - t^k) \geq 0;
    \end{equation*}
    \item For every $k$, the inequality $f(W^{k-1}) - f(W^k) - \|t^{k-1} - t^k\|^2\geq 0$ is valid;
    \item For each $k$, the inequality $f(W^*) \leq f(W^k) \leq f(W^*) + \|t^* - t^{k-1}\|^2$ holds true.
\end{enumerate}
\end{lemma}

\begin{proof}
The first statement follows the third statement of Lemma \ref{pro32} and the $\operatorname{FSLP}$ algorithm. The lower bound constraint \eqref{lower} provides a comprehensive way to prove this statement. The linearized bi-variant function at $t^{k-1}$ can be expressed as follows.
\begin{equation*}
\begin{split}
& t^\top CW(t)\mid_{t=t^{k-1}} \ge t^{{k-1}^\top}CW(t^{k-1}) + (t-t^{k-1})^\top CW(t^{k-1}) \\
& \qquad \qquad \qquad \qquad + t^{{k-1}^\top}C(W(t)-W(t^{k-1})) \\
\end{split}
\end{equation*}
Hence, at $t=t^k$, we have, 
\begin{equation*}
\begin{split}
& t^\top CW(t^k) \ge t^{{k-1}^\top}CW(t^{k-1}) + (t^k-t^{k-1})^\top CW(t^{k-1}) \\
& \qquad \qquad \qquad \qquad + t^{{k-1}^\top}C(W(t^k)-W(t^{k-1})) \\
\end{split}
\end{equation*}
It implies that, 
\begin{equation*}
\begin{split}
& (t^k-t^{k-1})^\top CW(t^k) \ge (t^k-t^{k-1})^\top CW(t^{k-1}) \\
& (t^k-t^{k-1})^\top (CW(t^k)-CW(t^{k-1})) \ge 0 \\
\end{split}
\end{equation*}
Finally, 
\begin{equation*}
\begin{split}
& (t^k-t^{k-1})^\top (t^{k+1}-t^k) \ge 0 \\
\end{split}
\end{equation*}

Now, let's examine the second statement. Since \(W^k\) is the optimal solution of \eqref{bi} for \(t = t^{k-1}\) and \(W^{k-1} \in \mathcal{F}\), it follows that:

\begin{equation*}
f(W^k)+|CW^k-t^{k-1}|^2 - f(W^{k-1})-|CW^{k-1}-t^{k-1}|^2\leq 0,
\end{equation*}
Factoring in that $CW^{k-1}=t^{k-1}$, we derive:
\begin{equation*}
f(W^{k-1})-f(W^k) - |t^k-t^{k-1}|^2\geq 0,\ \forall k\geq 1.
\end{equation*}

Pertaining to the third statement, given that \(W^*\) represents the global optimum for problem \eqref{newpro} and \(W^k\) is an element of \(\mathcal{F}\), we deduce the subsequent:

We can infer the following from \eqref{p31}:
\begin{equation}\label{p31}
f(W^k)- f(W^*)\geq 0, \quad \forall k.
\end{equation}
Also, given that $W^k$ is the optimal solution of problem \eqref{bi-new} when $t=t^k$ and $W^* \in \mathcal{F}$, we can deduce:
\begin{equation}\label{p32}
f(W^k)+|CW^k-t^{k-1}|^2\ - f(W^*) - \|CW^*-t^{k-1}\|^2 \leq 0.
\end{equation}
It then follows from the above two equations \eqref{p31} and \eqref{p32} that,
\begin{equation*}
f(W^*)\leq f(W^k)\leq f(W^*)+\|CW^*-t^{k-1}\|^2.
\end{equation*}
This concludes the proof.
\end{proof}

Subsequently, we will describe the characteristics of the accumulation point of a sequence that is produced by the $\operatorname{FSLP}$ algorithm. To achieve this, we present the concept of an $\epsilon$-local minimizer as defined below.

\begin{figure}[t]
    \centering
    \includegraphics[width= 0.9\columnwidth]{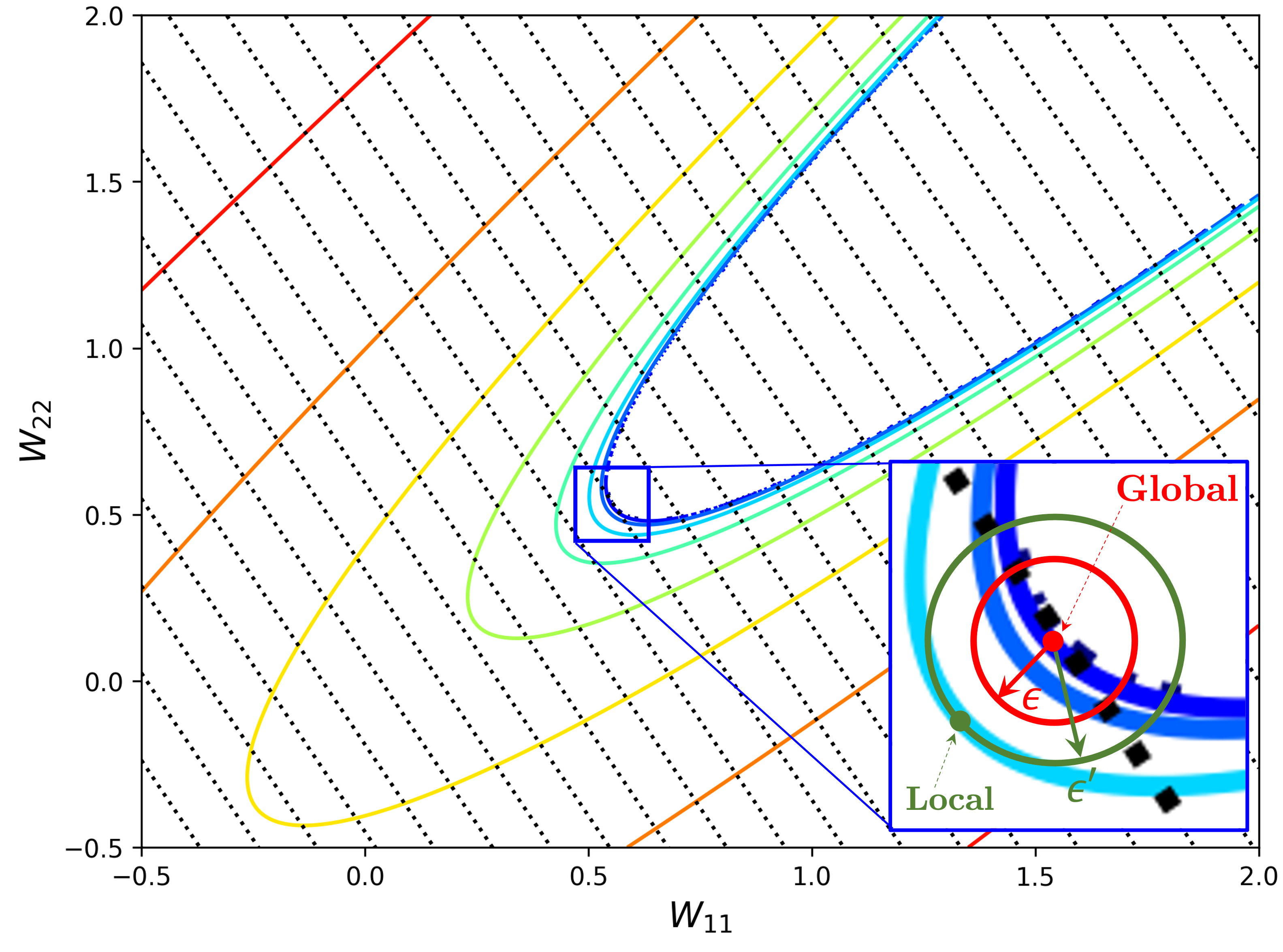}
    \caption{$\mathbf{\epsilon}$-approximation of global and local solutions. The green point is positioned outside the $\epsilon$ radius. The point is positioned in an expanded region with radius $\epsilon' > \epsilon$. The red point, representing the global solution, is located within the area bounded by the radius $\epsilon$.}
    \label{fig: fig6}
\end{figure}

\begin{definition}\label{defiii1}
We denote \(W^*\) as an \(\epsilon\)-local minimizer of problem \eqref{newpro} within a given neighborhood, given that \(W^*\) belongs to \(\mathcal{F}\),
\begin{equation*}
\mathcal{N}(W^*,\delta)={W:\|W-W^*\| \leq \sqrt{\epsilon}},
\end{equation*}
should the following condition hold for every \(W \in \mathcal{N}(W^*,\epsilon) \cap \mathcal{F}\):
\begin{equation}
f(W^*) - f(W)\leq \mathcal{O}(\epsilon)
\end{equation}
is satisfied, where $\lim _{\epsilon \rightarrow 0} \frac{\mathcal{O}(\epsilon)}{\epsilon}=\delta$ for certain value of $\delta>0$.
\end{definition}
Let's consider the following example that illustrates the concept of $\epsilon$-approximation in the vicinity of a global optimal solution.
\begin{equation}\label{example}
\begin{split}
\min_{W_{22},W_{11}} & \quad 100W_{11}^2+ 130W_{22}^2-210W_{11}W_{22}-60W_{11}+30 \\
s.t. &  \quad  W_{22}+2W_{11} \leq b_i, \quad \forall i \in \{1, 2, ..., 31\}\\
& \quad -0.5 \leq W_{ii} \leq 9.5, \quad \forall i \in \{1,2\} \\
& \quad W \in \mathcal{F}. 
\end{split}
\end{equation}
Where, $b_i$ is a linearspace on $[-1.5, 1.5]$, such as $b_i = -1.5 + \frac{i \times (1.5 - (-1.5))}{30}$, for $\forall i = \{1, 2, ..., 30\}$. Fig.~\ref{fig: fig6} compares the situation of two global and local solutions. 
In LCQP optimization problems, finding the exact optimal solution is often challenging, especially when the problem has many negative eigenvalues. Therefore, in practice, we often aim to find an approximate solution that is``good enough'' – a solution within a certain small distance, often denoted $\epsilon$, from the true global optimal solution. An $\epsilon$-approximation of global and local solutions refers to approximate solutions that lie within a distance $\epsilon$ from the true global and local optimal solutions.

Fig. \ref{fig: fig6} describes the points that represent solutions and circles with a certain radius to denote the $\epsilon$ neighborhoods (the areas within a distance $\epsilon$ from the points) around these solutions. The green point, which lies outside the $\epsilon$-radius, represents an approximation that does not meet the required level of precision. This point is situated in an expanded region with a radius of $\epsilon'$ (where $\epsilon' > \epsilon$), indicating that this approximation is further away from the optimal solution than we would ideally want.

However, the red point within the $\epsilon$-radius represents a more precise approximation. Since it is positioned within the area bounded by the $\epsilon$-radius around the true optimal solution, this indicates that the red point is close enough to the optimal solution to be considered a satisfactory $\epsilon$-approximation. It is important to note that the smaller the value of $\epsilon$, the closer the approximation is to the optimal solution. In other words, a smaller $\epsilon$ corresponds to a higher level of precision. Therefore, in this case, the red point represents a better approximation of the optimal solution than the green point.

The aforementioned definition can be regarded as an extensive interpretation of the term \(\epsilon\)-solution, as delineated in \cite{40}. As a result, we proceed to formulate the ensuing theorem.
\vspace{-0mm}

\begin{theorem} \label{t31}
Consider the case where \(W(t)\) represents the optimal solution to equation \eqref{bi} for a specific \(t \in [t_l, t_u]\). When the condition \(t = CW(t)\) is met, it follows that \(W(t)\) will act as an \(\epsilon\)-local solution to the function \(f(W)\) within a suitable neighborhood.
\begin{proof}
Given that \(W(t)\) is the optimal solution to problem \eqref{bi} and assuming the condition \(t = CW(t)\) is true, we can deduce:
\begin{equation}\label{proof1}
    f(W)+\|CW(t)-CW\|^2 \geq f(W(t)), \ \forall W\in \mathcal{F}.
\end{equation}
By leveraging Definition \ref{defiii1}, we can derive the following conclusion when a specific $\epsilon$ is set: 
\begin{equation*}
\begin{split}
    \|CW-CW(t)\|^2\leq \|C\|^2\|W-W(t)\|^2\leq \epsilon\|C\|^2, \\
\end{split}
\end{equation*}
When combined with \eqref{proof1}, it results in:
\begin{equation*}
    f(W)\geq f(W(t))-\epsilon.
\end{equation*}
where $\lim_{\epsilon \rightarrow 0} \frac{\mathcal{O}(\epsilon)}{\epsilon}=\|C\|^2$. This substantiates that $W(t)$ functions as a $\epsilon$-local solution to the problem \eqref{eq: 14}-\eqref{eq: 17}, with the neighborhood $\mathcal{N}(W(t), \epsilon)$ taken into account. This brings the theorem to the end of the proof.
\end{proof}
\end{theorem}

By employing Theorem \ref{t31} and Lemma \ref{le31}, we are equipped to deduce the convergence trajectories of the FSLP algorithm, aiming towards an \(\epsilon\)-local minimizer for problem \eqref{newpro}, delineated in the subsequent section:

\begin{theorem}
Suppose we possess a sequence \({(W^k, t^k)}\) formulated by the FSLP algorithm with \(\epsilon = 0\), and \((\bar{W}, \bar{t})\) serves as a convergence point for this sequence. In such a case, \(\bar{W}\) stands as not merely a KKT point but also an \(\epsilon\)-local minimizer for the given problem.
\begin{proof}
The second claim in Lemma \ref{le31} suggests that the sequence ${f(W^k)}$ is non-ascending. Given that $\mathcal{F}$ is a compact set and ${W^k}$ is a subset of $\mathcal{F}$, this means that the sequence ${f(W^k)}$ is bounded and therefore converges, leading to the situation where $f(W^{k-1})-f(W^k)$ approaches zero. The resulting conclusion is as follows:
\begin{equation}\label{proof22}
    \lim_{k\to\infty}\|t^{k-1}-t^k\|=0.
\end{equation}
Assume that $(\Bar{W},\bar{t})$ is an accumulation point of the sequence ${(W^k,t^k)}$. In such a case, there is a subsequence ${(W^{k_j},t^{k_j})}$ satisfying the following conditions:
\begin{equation*}
    W^{k_j}\rightarrow \Bar{W}, \ t^{k_j}\rightarrow \Bar{t}, \ where \ j\rightarrow \infty
\end{equation*}
From equation \eqref{proof22}, we can conclude that as $j$ approaches infinity, $t^{{k_j}-1}$ approaches $\Bar{t}$. Given that $t^k=CW^k$ holds for all $k$, it follows that $\bar t=C\bar W$. The compactness of $\mathcal{F}$ implies that $\bar W$ is in $\mathcal{F}$. Given that $W^{k_j}$ is an optimal solution to problem \eqref{bi} with $t=t^{{k_j}-1}$, we can deduce:
\begin{multline*}
     W^\top A_+W+b^\top W-2(t^{{k_j}-1})^\top CW+\|t^{{k_j}-1}\|^2\\
      \geq (W^{k_j})^\top A_+W^{k_j}+b^\top W^{k_j}-2(t^{{k_j}-1})^\top CW^{k_j}+\|t^{{k_j}-1}\|^2, \\
     \forall W \in \mathcal{F}.
\end{multline*}
By applying the limit to both sides of the inequality as $j$ tends towards infinity ($j\rightarrow \infty$), we can derive:
\begin{align*}
     & \ W^\top A_+W+b^\top W-2\bar{t}^\top CW+\|\bar{t}\|^2\\
     & \geq \bar{W}^\top A_+\bar{W}+b^\top \bar{W}-2\bar{t}^\top C\bar{W}+\|\bar{t}\|^2, \ \forall W \in \mathcal{F}.
\end{align*}
This implies that \(\bar{W}\) serves as the optimal solution for equation \eqref{bi} where \(t=\bar{t}=C\bar{W}\), and concurrently acts as an \(\epsilon\)-local solution as well as a Karush-Kuhn-Tucker (KKT) point for the problem in question. Hence, the proof reaches its conclusion.
\end{proof}
\end{theorem}

Nonetheless, in scenarios where \(r > 1\), the traditional Successive Linear Programming (SLP) algorithm may face challenges in pinpointing the optimal solutions for problem \eqref{newpro}. This difficulty stems from the fact that the conventional SLP algorithm depends on scrutinizing the convex hull formed by the accumulation points generated by the algorithm itself. Furthermore, the bounded $\operatorname{FSLP}$ algorithm can pinpoint a highly satisfactory feasible solution, which can be employed as an upper limit in the novel Branch-and-Bound (B\&B) method.

\section{A novel branch-and-bound Algorithm}

In this part, we propose a new comprehensive algorithm designed to find a globally optimal solution for problem \eqref{eq:14}, while maintaining compliance with a predetermined \(\epsilon\)-tolerance. This objective is accomplished by amalgamating the \(\operatorname{FSLP}\) algorithm, the branch-and-bound (B\&B) structure, and convex relaxation techniques. Additionally, we guarantee the convergence of our proposed algorithm and examine its complexity.

\subsection{Tighter Lifted Convex Relaxation}
Convex relaxation of the lifted problem is a methodology employed in optimization theory. This technique involves simplifying a complex or nonconvex problem---referred to as the "lifted" problem, which signifies a problem that has been altered or enhanced to facilitate easier resolution---by approximating it to a more manageable convex problem. This simplification aids in the problem's resolution, typically utilizing approaches that ensure the identification of the global optimum, thereby preventing the solution process from getting entrapped in local optima.

We initiate by looking at a restricted adaptation of the lifted problem, where the variable $t$ is bound by a subbox $[l, u]$, with both $l$ and $u$ belonging to $\mathbb{R}^r$.
\begin{equation}\label{res}
\begin{split}
& \min \quad  f(W,t)=W^\top A_{+}W+b^\top W-\|t\|^2 \\
& \ s.t. \quad \ CW-t=0,\\
& \quad \qquad W \in \mathcal{F}. \qquad t\in [l, u]
\end{split}
\end{equation}
where $l, u \in \mathbb{R}^r$ with $[l, u] \subseteq [t_l, t_u]$.

Denote the $i^{th}$ row vector of the matrix $C$ as $c_{i}^\top$, with $\lambda_i$ representing the absolute value of the negative eigenvalues of $A$ for all $i=\{1,...,r\}$, and $p_i$ as the respective orthogonal unit eigenvectors for all $i=\{1,...,r\}$. According to the singular value decomposition of $A$, we get $t_i=c_{i}^\top W=\sqrt{\lambda_i}p_{i}^\top W$ for all $i=\{1,...,r\}$. Given that $\begin{matrix} \sum_{i=1}^n p_{i}p_{i}^\top=I \end{matrix}$, we can deduce:
\begin{equation*}
\begin{split}
 \sum_{i=1}^r \frac{s_i}{\lambda_i} & = \sum_{i=1}^r W^\top \left (p_{i}p_{i}^\top \right )W \leq W^\top \left (\sum_{i=1}^r p_{i}p_{i}^\top \right )W \\
& =W^TW \leq \Bar{u}^\top W, \qquad \forall W\in \mathcal{F}\
\end{split}
\end{equation*}
where $s_i$ equals the square of $t_i$, and $\Bar{u}$ is the upper boundary of the variable $x$ over the set $\mathcal{F}$. Meanwhile, when $t$ falls within the range $[l,u]$, we have:
\begin{equation*}
    s_i=t_{i}^2 \leq (l_i+u_i)t_i-l_{i}u_{i}, \qquad i=\{1,...,r\}. 
\end{equation*}
Thus, we can infer the subsequent convex relaxation for the equation \eqref{res}:
\begin{equation} \label{cr}
\begin{split}
& \min \quad  f(W,t)=W^\top A_{+}W+b^TW-\sum_{i=1}^r s_i, \\
& \ s.t. \quad \ CW-t=0,\\
& \quad \qquad t_{i}^2 \leq s_i, \\
& \quad \qquad s_i \leq (l_i+u_i)t_{i}-l_{i}u_{i}, \\
& \quad \qquad \sum_{i=1}^r \frac{s_i}{\lambda_i} \leq \Bar{u}^\top W, i=1,...,r,\\
& \quad \qquad W \in \mathcal{F}, t\in [l, u].
\end{split}
\end{equation}

Following this, we examine the objective values at the optimal solutions for problem \eqref{res} and its relaxation \eqref{cr}. It is found that:
\begin{theorem} \label{t51}
Consider $f_{[l,u]}^*$ and $v^*_{[l,u]}$ as the optimal values of \eqref{res} and its relaxed counterpart \eqref{cr}, respectively. Let's assume $(\bar{W},\bar{s},\bar{t})$ to be the optimal solution to \eqref{cr}. In that case, we get:
\begin{equation*}
    0\leq f_{[l,u]}^*-v_{[l,u]}^*\leq f(\bar{W},\bar{t})-v^*_{[l,u]}\leq \frac{{\|u-l\|}^2}{4}.
\end{equation*}
\begin{proof}
Clearly, $(\bar{W},\bar{t})$ is a feasible solution to \eqref{res}. Since \eqref{cr} is a convex relaxation of \eqref{res}, it follows from the choice of $f_{[l,u]}^*$ that
\begin{align}\label{proof3}
    0 & \leq f_{[l,u]}^*-v^*_{[l,u]}\leq f(\bar{W},\bar{t})-v^*_{[l,u]} \nonumber \\
      & = \sum_{i=1}^r(\bar{s}_i-\bar{t}_i^2) \\
      & \leq \sum_{i=1}^r[-\bar{t}_i^2+(l_i+u_i)\bar{t}_i-l_iu_i] \nonumber \\
      & \leq \frac{1}{4}\sum_{i=1}^r(u_i-l_i)^2=\frac{1}{4}\|u-l\|^2,\nonumber
\end{align}
with the latter inequality being a direct consequence of the constraints placed upon \(s_i\) in equation \eqref{cr}. This brings us to the conclusion of our proof.
\end{proof}
\end{theorem}

In light of the demonstration provided for Theorem \ref{t51}, we can propose the following:

\begin{proposition} \label{p51}
Given that \((\bar{W}, \bar{s}, \bar{t})\) represents the optimal solution for equation \eqref{cr} and \(\epsilon > 0\). If the condition \(\sum_{i=1}^r [\bar{s}_i - \bar{t}_i^2] \leq \epsilon\) is satisfied, then the pair \((\bar{W}, \bar{t})\) qualifies as an \(\epsilon\)-optimal solution to equation \eqref{res}.
\end{proposition}

\begin{proof}
We observe that \((\bar{W}, \bar{t})\) is a feasible solution to problem \eqref{res}. From equation \eqref{proof3}, we have
\begin{equation*}
    0 \leq f(\bar{W}, \bar{t})-f_{[l,u]}^* \leq f(\bar{W}, \bar{t})-v^*_{[l,u]} \leq \sum_{i=1}^r (\bar{s}_i-\bar{t}_i^2) \leq \epsilon.
\end{equation*}
Therefore, we can deduce
\begin{equation*}
    f(\bar{W}, \bar{t}) \leq f_{[l,u]}^* + \epsilon,
\end{equation*}
signifying that \((\bar{W}, \bar{t})\) qualifies as an \(\epsilon\)-optimal solution to problem \eqref{res}.
This brings our proof to a close.
\end{proof}

Observe that, given a feasible point \(W^*\), we have the ability to meticulously select a neighborhood \(\mathcal{N}(W^*, \epsilon)\) that satisfies the subsequent constraint,
\begin{equation} \label{cons}
    t_i \in \left [l_i=(CW^*)_i-\sqrt{\epsilon/r},\ u_i=(CW^*)_i+\sqrt{\epsilon/r} \right ],
\end{equation}
In the end, the bound difference can be characterized as a function of the number of eigenvalues and the tolerance level established for the approximation in the quest for the global solution.  
\begin{equation}\label{bound}
u_i-l_i - 2 \sqrt{\epsilon / r }\leq 0, \quad i=\{1, \ldots, r\}
\end{equation}
From Theorem \ref{t51}, we can obtain the following Theorem.

\begin{theorem} \label{t52}
Consider \(W^* \in \mathcal{F}\) and a neighborhood \(\mathcal{N}(W^*, \delta)\) that fulfills all the constraints delineated in \eqref{cons}. Let \((\bar{W}, \bar{s}, \bar{t})\) denote the optimal solution to the relaxed problem \eqref{cr}. In this case, \((\bar{W}, \bar{t})\) stands as an \(\epsilon\)-local solution to problem \eqref{res}, relative to the neighborhood \(\mathcal{N}(W^*, \epsilon)\).
\end{theorem}

\begin{proof}
Let us denote the optimal values of \eqref{res} and its relaxed counterpart \eqref{cr} as \(f_{[l,u]}^*\) and \(v^*_{[l,u]}\), respectively. Given Theorem \ref{t51}, we can infer that
\begin{equation*}
    f(\bar{W}, \bar{t}) - f_{[l,u]}^* \leq f(\bar{W}, \bar{t}) - v^*_{[l,u]} \leq \frac{\|u - l\|^2}{4}.
\end{equation*}

Together with \eqref{cons}, we derive that
\begin{equation*}
    f(\bar{W}, \bar{t}) - f_{[l,u]}^*\leq \frac{\|u - l\|^2}{4} \leq \min_{W \in \mathcal{N}(W^*, \epsilon), t \in [l,u]} f(W, t) + \epsilon.
\end{equation*}
Hence, \((\bar{W}, \bar{t})\) qualifies as an \(\epsilon\)-local solution to \eqref{res} in relation to the neighborhood \(\mathcal{N}(W^*, \epsilon)\). This brings the proof to its conclusion.
\end{proof}

Theorem \ref{t52} sets the foundation for segmenting each interval $[t_{l}^i, t_{u}^i]$ into discrete subintervals, of which the count can be determined by $\big \lceil \frac{\sqrt{r}(t_{u}^i-t_{l}^i)}{2\sqrt{\epsilon}} \big \rceil$. These intervals are represented as $[l_i, u_i]$. The definition of $t_{l}^i$ and $t_{u}^i$ are guided by their respective formulae.
This process results in the decomposition of the feasible domain, $\mathcal{F}$, into a multitude of smaller, manageable subdomains. Each one of these subdivisions is characterized by its unique box constraints, which impose restrictions such that $t_i$ lies within $[l_i, u_i]$ and the constraint width, $u_i-l_i$, does not exceed $2\sqrt{\epsilon /r}$.

After implementing this subdivision, we can address problem \eqref{cr} within each of these distinct areas, thereby acquiring an \(\epsilon\)-local approximation for every individual region. Subsequently, from all these locally approximated solutions, we select the one that results in the smallest objective function value to serve as our final solution.

Consequently, the final solution obtained inherently qualifies as a global \(\epsilon\)-approximate solution to the initial problem \eqref{newpro}. Therefore, this segmentation strategy effectively furnishes a viable \(\epsilon\)-approximation for \eqref{newpro}.

For ease of reference, we term this procedure as the ``brutal force algorithm''. The ``Brute force'' refers to a general problem-solving approach in computer science and operations research. A brute force algorithm is a simple, straightforward approach to solve a problem based on the problem's definition and feasible solutions. It generates all possible solutions for a problem and checks each solution individually to see if it satisfies the problem's constraints. If the solution satisfies the constraints, it is considered a valid solution; if not, it is discarded, and the algorithm moves on to the next potential solution. The forthcoming theorem validates the aforementioned assertions.
\begin{theorem}
The "brutal force algorithm" is capable of discovering a global $\epsilon$-approximate solution to \eqref{newpro} within a time complexity of $\mathcal{O}\left (\prod_{i=1}^{r} \big \lceil \frac{\sqrt{r}(t_{u}^i-t_{l}^i)}{2\sqrt{\epsilon}} \big \rceil N \right)$, where $N$ represents the complexity involved in the resolution of \eqref{cr}.
\end{theorem}

While providing an $\epsilon$-approximate solution, the brutal force algorithm is generally inefficient due to its conservative approach to working with small neighborhoods. Essentially, the algorithm partitions the problem into many tiny subproblems and tries to find local solutions within each one, which could lead to a high computational burden, especially for larger problems. It can be deemed an overkill approach due to the excessive fine-grain partitioning involved.

However, it provides a solid theoretical foundation for building more efficient algorithms. In the following subsection, we will explore how to improve upon this by integrating the Faster $\operatorname{FSLP}$ algorithm with other methods. The goal is to develop an enhanced global algorithm for the problem expressed in \eqref{newpro}. This approach is expected to be more effective and less computationally demanding and will allow us to handle more complex problems in a reasonable amount of time.

\subsection{The $\operatorname{FSLP}$-embedded Branch-and-Bound algorithm}
In this section, we delve into the architecture of a newly formulated global algorithm that adeptly integrates the Faster Successive Linear Programming ($\operatorname{FSLP}$) approach principles with the well-established branch-and-bound (B\&B) techniques. The application of convex quadratic relaxation strategies enhances this integration. For simplicity, we label this innovative hybrid algorithm as $\operatorname{FBB}$.

$\operatorname{FBB}$ is endowed with two key characteristics that empower it with a distinctive computational edge:
\begin{enumerate}
    \item \emph{Global Optimality Verification and Upper Bound Refinement:} The $\operatorname{FBB}$ algorithm can verify the global optimality of the solutions procured through $\operatorname{FSLP}$ and improve the existing upper bound by invoking a new $\operatorname{FSLP}$ process under specific conditions. This unique feature allows greater flexibility and precision in searching for the global optimum.
    \item \emph{Acceleration of B\&B Convergence Using the $\operatorname{FSLP}$ Upper Bound:} $\operatorname{FBB}$ leverages the upper bound derived from the $\operatorname{FSLP}$ method to accelerate the convergence rate of the B\&B approach. This feature significantly reduces the computational resources required, making $\operatorname{FBB}$ a much more efficient tool to tackle the problem \eqref{newpro}.
\end{enumerate}

Using the strength of these characteristics, $\operatorname{FBB}$ provides a robust and effective approach to solving complex optimization problems such as \eqref{newpro}. In the subsequent sections, we will outline the specifics of these features and their application within the $\operatorname{FBB}$ framework.

Let's begin by exploring the novel branching technique, which is fundamentally predicated on dividing the feasible rectangle (expressed as $\Delta$) of $t$ in equation \eqref{cr}. The domain $\Delta$ is split into multiple subrectangles during each iteration. A conventional approach is to partition the longest side of the subrectangle at its center point. However, this method might not effectively enhance the lower bound as it does not exploit the particular structure of equation \eqref{cr}.

To address this limitation, we propose an adaptive branch-and-cut rule. Take into consideration a sub-rectangle $\Delta$ represented as $[l,u]$, and let the optimal solution to equation \eqref{cr} over $\Delta$ be $(W^*, s^*, t^*)$. If condition $\sum_{i=1}^r s_i^* - (t_i^*)^2 \leq \epsilon$ holds, then by Proposition, $W^*$ is considered an $\epsilon$-optimal solution to \eqref{res} over $\Delta$, and further partitioning of this subrectangle will be unnecessary.

On the other hand, if the aforementioned condition is not met, it signifies the existence of at least one \(i\) within the set \(\{1, \ldots, r\}\) such that \(s_i^* - (t_i^*)^2 \geq \frac{\epsilon}{r}\). In this case, we select \(i^*\) which minimizes the expression \(s_i^* - (t_i^*)^2\). Subsequently, the sub-rectangle \(\Delta\) is divided into two new sub-rectangles, \(\Delta_1\) and \(\Delta_2\), along the edge \([l_{i^*}, u_{i^*}]\) at a specific point \(w\). This point \(w\), lying in the interval \((l_{i^*}, u_{i^*})\), is referred to as the branching value. Following this, the constraints for the intervals \([l_{i^*}, w]\) and \([w, u_{i^*}]\) are considered individually.
\begin{equation}\label{cons}
    s_{i^*}\leq (l_{i^*}+w)t_{i^*}-l_{i^*}w, \quad
    s_{i^*}\leq (w+u_{i^*})t_{i^*}-wu_{i^*}
\end{equation}

\begin{figure}[t]
\centerline{\includegraphics[width=0.7\textwidth]{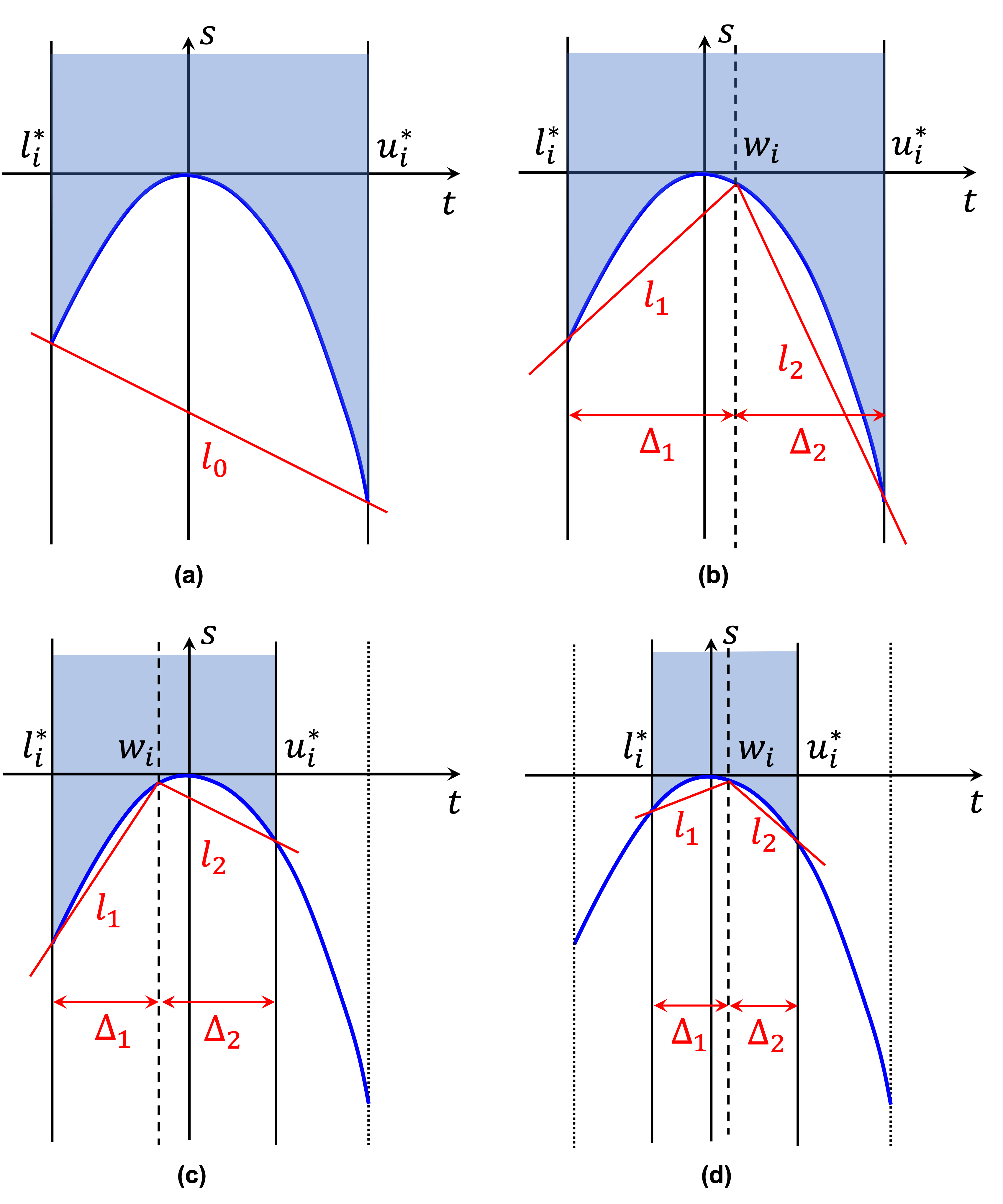}}
\caption{Partition on the variable $t_{i^*}$.}
\label{fig: fig7}
\end{figure}

Referring to Fig.~\ref{fig: fig7}, it is clear that the two constraints from equation \eqref{cons} effectively eliminate $(W^*, s^*, t^*)$ when $w$ is equivalent to $t_{i^*}^*$. Let's define two linear equations, $l_1$ and $l_2$, given by $l_1: s_{i^*}=(l_{i^*}+w)t_{i^*}-l_{i^*}w$ and $l_2: s_{i^*}=(w+u_{i^*})t_{i^*}-wu_{i^*}$ respectively. To identify the optimal linear cuts, we formulate an optimization problem (referenced as \eqref{eq: area}). The goal of this problem is to minimize the area confined by the parabolic curve $s_i^*=(t_i^*)^2$ and the two straight lines $l_1$ and $l_2$, within the interval $[l_{i^*}, u_{i^*}]$. 

\begin{center}
\hspace*{-1cm} 
\begin{minipage}{1.15\textwidth} 
\begin{algorithm}[H]
\caption{The $\operatorname{FBB}$ Algorithm}
\label{a2}

\KwIn{Initialized with inputs $A, b, L, l$ and a stop criterion $\epsilon > 0$}
\KwOut{A global $\epsilon$-solution denoted by $W^*$}

\textbf{Step 0 (Preparation)}: Identify $r$, the eigenvalues of $A$. If $r > 5$, set $\rho = 2^r$ and for each $j \in \{1, \ldots, \rho\}$, let $\mu_j$ be a binary vector of size $r$ with elements in $\{-1,1\}$. If $r > 5$, limit $\rho$ to $2$. Consider $\mu_1$ and $\mu_2$ as $r$-dimensional vectors, where every element of $\mu_1$ is $1$, and every element of $\mu_2$ is $-1$. For each $j \in \{1, \ldots, \rho\}$, solve the problem $\min_{(W,t) \in \mathcal{F}_t} \mu_j^T t$ to derive the optimal solution, denoted as $t_j$\;

\textbf{Step 1 (Local Solutions)}: Identify nearly local $\epsilon$-optimal solutions, referred to as $W_j^k$, for the problem \eqref{newpro} by executing the $\operatorname{FSLP}\left(t^0, \epsilon\right)$ function, where $t^0 = t_j$ for each $j \in \{1, \ldots, \rho\}$. Determine the argument $W^*$ that yields the minimum function value of the set $\{f(W_j^k)~|~j=1, \ldots, \rho\}$. The value of the function in this argument, $f(W^*)$, is assigned to $v^*$\;

\textbf{Step 2 (Optimal Solution and Initialization)}: Solve the problem \eqref{cr} within the boundaries $[l=t_l, u=t_u]$ to yield an optimal solution $(W^0, s^0, t^0)$, along with a lower limit $v^0$. If $f(W^0)<v^*$, update the upper limit $v^*$ with $f(W^0)$ and replace $W^*$ with $W^0$. Set $k:=0$, define $l^k$ and $u^k$ as $t_l$ and $t_u$ respectively. Set the region $\Delta_k$ as the interval $[l^k, u^k]$. Initialize the set $\Omega$ as $\{[\Delta_k, v^k, (W^k, s^k, t^k)]\}$\;

\textbf{Step 3}: \While{$\Omega \neq \emptyset$}{

\textbf{Phase (I) Selection Node}: From the set $\Omega$, select the node $\left[\Delta_k, v^k,\left(W^k, s^k, t^k\right)\right]$ with the smallest value of $v^k$. Remove the node from $\Omega$\;

\textbf{Phase (II) Termination}: If $v^k \geq v^* -\epsilon$ or $f\left(W^k\right) - v^k \leq \epsilon$, then both $W^*$ and $W^k$ qualify as $\epsilon$ optimal solutions to \eqref{newpro}, indicating the algorithm ends\;

\textbf{Phase (III) Partitioning}: Identify the index $i^*$ that maximizes $s_i^k - (t_i^k)^2$ for every $i \in \{1, \ldots, r\}$. Set $w_{i^*}$ as the midpoint of the interval $[l_{i^*}^k, u_{i^*}^k]$, which equals $\frac{l_{i^*}^k + u_{i^*}^k}{2}$. If the pair $\left(s_{i^*}^k, t_{i^*}^k\right)$ is an element of the set $\Gamma_k\left(w_{i^*}\right)$, assign the branching point $\beta_{i^*}$ to $w_{i^*}$. Otherwise, let $\beta_{i^*}$ be $t_{i^*}^k$. Divide the rectangle $\Delta_k$ into two smaller rectangles, $\Delta_k^{k_1}$ and $\Delta_k^{k_2}$, along the edge $\left[l_{i^*}^k, u_{i^*}^k\right]$ at the branching point $\beta_{i^*}$\;

\textbf{Phase (IV) Problem Solving and Set Expansion}: Solve the problem \eqref{cr} over each of the two new subrectangles $\Delta_k^{k_1}$ and $\Delta_k^{k_2}$. This yields a lower bound, denoted as $v^{k_j}$, as well as an optimal solution $\left(W^{k_j}, s^{k_j}, t^{k_j}\right)$ for each $j=\{1,2\}$. Append to the set $\Omega$ the tuples $\left[\Delta_k^{k_1}, v^{k_1}, \left(W^{k_1}, s^{k_1}, t^{k_1}\right)\right]$ and $\left[\Delta_k^{k_2}, v^{k_2}, \left(W^{k_2}, s^{k_2}, t^{k_2}\right)\right]$\;
}
\end{algorithm}
\end{minipage}
\end{center}

\begin{algorithm}[H]
\caption{Algorithm Description}
\While{true}{
    new subrectangles, $\Delta_k^{k_1}$ and $\Delta_k^{k_2}$. Doing this yields a lower bound, denoted as $v^{k_j}$, as well as an optimal solution that is represented as $\left(W^{k_j}, s^{k_j}, t^{k_j}\right)$ for each $j=\{1,2\}$. After obtaining these values, expand the set $\Omega$ by including the results related to the new sub-rectangles. Specifically, append to the set $\Omega$ the tuples $\left[\Delta_k^{k_1}, v^{k_1}, \left(W^{k_1}, s^{k_1}, t^{k_1}\right)\right]$ and $\left[\Delta_k^{k_2}, v^{k_2}, \left(W^{k_2}, s^{k_2}, t^{k_2}\right)\right]$, which contain the information about the newly created subrectangles, their corresponding lower bounds, and optimal solutions.\\
    \textbf{Phase (V) (Reinitiating $\operatorname{FSLP}$ and Updating the Solution)}: First, define $\hat{W}$ as the solution that minimizes the function $f$ over the optimal solutions of the new subrectangles, that is, $\hat{W} = \arg \min \{f(W^{k_1}), f(W^{k_2})\}$.\\
    \If{$f(\hat{W}) < v^*$}{
        This implies that the function value at $\hat{W}$ is less than our current best upper bound $v^*$. Under this condition, execute the function $\operatorname{FSLP}(t^0, \epsilon)$ with an initial value $t^0 = C \hat{W}$ to find an $\epsilon$-local solution $\bar{W}^k$ for problem \eqref{newpro}.\\
        Following this, update the current optimal solution $W^*$ and the upper bound $v^*$. The updated $W^*$ becomes the argument that minimizes the function $f$ over $\hat{W}$ and $\bar{W}^k$, represented as $W^* = \arg \min \{f(\hat{W}), f(\bar{W}^k)\}$. The upper bound, denoted as $v^*$, is updated to reflect the value of the function evaluated at the new optimal solution, expressed as $v^* = f(W^*)$.
    }
    \textbf{Phase (VI) (Node Pruning)}: Here, we remove any redundant or non-competitive nodes from the set $\Omega$. Specifically, remove all nodes denoted by $\left[\Delta_k^j, v^j, \left(W^j, s^j, t^j\right)\right]$ for which $v^j \geq v^* - \epsilon$. Nodes that meet this criterion are less likely to lead to an improved solution. This is because their lower bound, denoted by $v^j$, is not significantly lower than the current best solution, $v^*$, given the specified tolerance, $\epsilon$. Following the pruning process, increment the iteration counter $k$ by one, which means moving on to the next iteration of the algorithm.
}
\end{algorithm}

Now, we are prepared to delve into the specifics of the innovative global algorithm designed to deal with problem \eqref{newpro}. The answer to this straightforward unconstrained optimization problem is given as follows:
\begin{align*}\
    & S(w) = \frac{1}{2} (u_{i^*}-l_{i^*})w^{2}-\frac{1}{2} (u_{i^*}^2-l_{i^*}^2)w+\frac{1}{6} (u_{i^*}^3-l_{i^*}^3).\\
\end{align*}
Figure~\ref{fig: fig7} illustrates that choosing the midpoint, $w^*=\frac{l_{i^*}+u_{i^*}}{2}$, as the point of partition along the edge $[l_{i^*}, u_{i^*}]$ of the sub-rectangle $\Delta$, minimizes the area $S(w)$. This strategy is preferred for dividing $\Delta$. However, if implementing the secant constraints at $w=w^*$, as shown in equation \eqref{cons}, does not eliminate the solution set $(W^*, s^*, t^*)$ for equation \eqref{cr} within $\Delta$, the lower bound improvements for either sub-rectangle $\Delta_1$ or $\Delta_2$ could be minimal.

To bypass this issue, an alternative approach is to select $w=t_{i^*}^*$. Specifically, if the constraints defined in equation \eqref{cons} with $w=w^*$ successfully exclude the solution set $(W^*, s^*, t^*)$, then $\Delta$ should indeed be partitioned at the midpoint $w^*$. Otherwise, to effectively refine the bounding process, partitioning $\Delta$ at $w^*=t_{i^*}^*$, based on the optimal value of the decision variable $t_{i^*}^*$, becomes the preferred strategy.



Now, we are prepared to delve into the specifics of the innovative global algorithm designed to deal with problem \eqref{newpro}.
In Algorithm \ref{a2}, Part III, \(\Gamma_k(w_{i^*})\) can be represented by two secant cuts:

$$
\Gamma_k\left(w_{i^*}\right)=\left\{\begin{array}{l|l}
(s_{i^*}, t_{i^*}) & \begin{array}{l}
s_{i^*}>(l_{i^*}^k+w_{i^*}) t_{i^*}-l_{i^*}^k w_{i^*} \\
s_{i^*}>(w_{i^*}+u_{i^*}^k) t_{i^*}-w_{i^*} u_{i^*}^k
\end{array}
\end{array}\right.
$$

The following observations pertain to the stages of the $\operatorname{FBB}$ algorithm:
\begin{enumerate}
    \item To begin, we employ the $\operatorname{FSLP}$ algorithm, utilizing various feasible points from the higher-level problem as starting points, with the aim of discovering a satisfactory \(\epsilon\)-local solution for equation \eqref{newpro}.
    \item During the first phase of the third step, the \(\operatorname{FSLP}\) algorithm is reinitialized to locate a more favorable \(\epsilon\)-local solution, on the condition that the objective value at the feasible point derived from the solution of the relaxation problem is lower than the current upper limit.
    \item Subsequently, in the third phase of the third step, we proceed to prune the optimal solution for the relaxation problem after each iteration in order to improve the lower bound.
\end{enumerate}
We will now demonstrate the convergence of the $\operatorname{FBB}$ Algorithm. Let's consider a sub-rectangle, denoted as $\bar{\Delta}=[\bar{l}, \bar{u}]$, that complies with the following condition.
\begin{equation}\label{satis}
    \bar{u}_i - \bar{l}_i - 2\sqrt{\epsilon /r} \leq 0 \qquad \forall i \in \{1,...,r\}.
\end{equation}
Assume that $(\bar{W}, \bar{s}, \bar{t})$ is the optimal solution for equation \eqref{cr} over the range of $\bar{\Delta}$. As stated in equation \eqref{proof3}, we find that,
\begin{equation}\label{note}
    \bar{s}_i-\bar{t}_{i}^2 \leq \frac{1}{4} (\bar{u}_i-\bar{l}_i)^2, i=\{1,...,r\},
\end{equation}
This suggests that,
\begin{equation} \label{cond}
    \sum_{i=1}^{r} [\bar{s}_i-(\bar{t}_i)^2] \leq \epsilon.
\end{equation}
Based on Proposition \ref{p51}, we can affirm that \(\bar{W}\) is an \(\epsilon\)-optimal solution for equation \eqref{res} within the domain of \(\bar{\Delta}\). In alignment with the third and fourth procedures of the \(\operatorname{FBB}\) algorithm, this specific sub-rectangle \(\bar{\Delta}\) will be conclusively removed from the set \(\Omega\). This signifies that a sub-rectangle $\Delta$ that complies with equation \eqref{satis} will no longer be partitioned in the course of the algorithm.

In the event that the optimal solution \((\bar{W}, \bar{s}, \bar{t})\) of a current sub-problem fails to meet the requirements of equation \eqref{cond} and necessitates further investigation, there must exist an index \(i^*\) within the set \(\{1, \ldots, r\}\) such that \(\bar{s}_{i^*} - \bar{t}_{i^*}^2 > \frac{\epsilon}{r}\). As per the equation, this implies that \(\bar{u}_{i^*} - \bar{l}_{i^*} > 2\sqrt{\frac{\epsilon}{r}}\), signifying that the length of the \(i^*\)-th edge of the corresponding sub-rectangle must be greater than \(2\sqrt{\frac{\epsilon}{r}}\).

As specified in Step 3, Phase III, it will proceed to division at the point $\bar{t}_{i^*}$ or at the midpoint of this edge. In addition, the length of any edge of the subrectangle associated with a subproblem that needs to be solved will always be greater than or equal to $2\sqrt{\epsilon /r}$.

Consequently, every edge of the starting rectangle is divided into no more than $\left\lceil \frac{\sqrt{r}(t_{u}^i-t_{l}^i)}{2\sqrt{\epsilon}} \right\rceil$ segments. This means that to obtain an $\epsilon$-accurate solution for equation \eqref{newpro}, the quantity of simplified subproblems that must be addressed with the $\operatorname{FBB}$ algorithm is bounded by,
\begin{equation*}
    \prod_{i=1}^{r} \big \lceil \frac{\sqrt{r}(t_{u}^i-t_{l}^i)}{2\sqrt{\epsilon}} \big \rceil.
\end{equation*}
Based on the prior discussion, it can be concluded that,
\begin{theorem}
The $\operatorname{FBB}$ Algorithm is capable of discovering a global $\epsilon$-approximate solution to the given problem within a maximum time frame of $\mathcal{O}(\prod_{i=1}^r [\frac{\sqrt{r}(t_u^i-t_l^i)}{2\sqrt{\epsilon}}]N)$, where $N$ signifies the complexity associated with solving problem \eqref{cr}.
\end{theorem}

\section{Computational Results}

The effectiveness of the $\operatorname{FSLP}$ and $\operatorname{FBB}$ algorithms for providing precise solutions for instances in the \texttt{PG-lib} library \cite{41} has been illustrated through numerous tests and a comparison made with the nonlinear solver \texttt{Baron 23.6.15} \cite{42} and the nonconvex mode of the \texttt{Gurobi 10} \cite{46} solver. These tests were carried out on a server equipped with dual Intel Xeon CPUs, each with 12 cores and 2 threads operating at $2.5~\mathrm{GHz}$, accompanied by $4 \times 16\mathrm{GB}$ RAM and a Linux-based operating system. A maximum computation time of $3$ hours was permitted for all methods. Baron utilized the multithreaded variant of \texttt{CPLEX20.10} \cite{43}, supporting up to 64 threads. On the other hand, the $\operatorname{FSLP}$ and $\operatorname{FBB}$ methods incorporated the \texttt{Mosek} \cite{44} semidefinite solver to solve Penalized semidefinite programs. At each branch and bound node, distinct solvers were used for solving the LCQP of the $\operatorname{FBB}$ method and the \texttt{CPLEX 20.10} solver was employed for solving the QP of the $\operatorname{FSLP}$ method. For computing feasible local solutions, \texttt{Ipopt} \cite{45} was utilized as a local solver.

The tests were conducted utilizing data from medium-sized power networks ranging from $3$ to $300$ buses, following the ACOPF formulation defined by constraints \eqref{eq: 14}-\eqref{eq: 16}. Table \ref{tab: Table1} offers comprehensive characteristics for each test case, including the name of the cases, the number of buses, generators, and lines in the examined power network. It also displays the optimal solution discerned by the $\operatorname{FSLP}$ and $\operatorname{FBB}$ methods within the allotted 3 hours of computing time.

\begin{figure}[t]
\centerline{\includegraphics[width=0.9\textwidth]{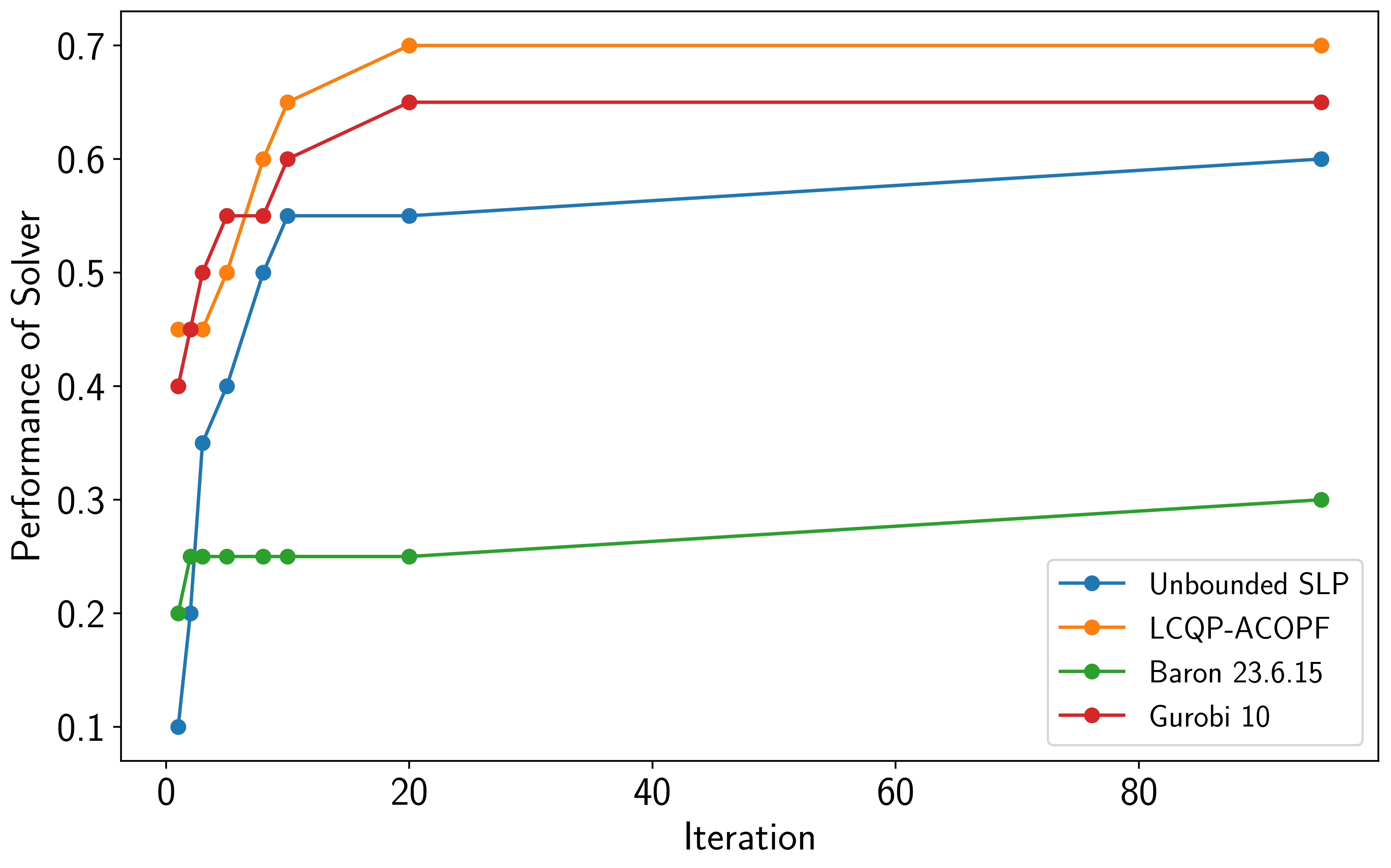}}
\caption{The solver performance profile concerning the cumulative duration for networks, ranging from 2 to 300 bus systems, within a time boundary of 3 hours for the proposed LCQP-ACOPF model, presented ACOPF model with unbounded-SLP in \cite{47}, Baron and Gurobi.}
\label{fig: fig11}
\end{figure}

\begin{table}[h]
\centering
\begin{tabular}{|l|c|c|c|c|}
\hline 
Name & $|\mathcal{N}|$ & $|\mathcal{G}|$ & $|\mathcal{L}|$ & Opt \\
\hline 
\texttt{caseWB2} & 2 & 1 & 1 & 12.262 \\
\hline 
\texttt{caseWB3} & 3 & 1 & 2 & 587.340 \\
\hline 
\texttt{pglib\_case3\_lmbd} & 3 & 3 & 3 & 5655.465 \\
\hline 
\texttt{caseWB5} & 5 & 2 & 6 & 9.076 \\
\hline 
\texttt{pglib\_case5\_pjm} & 5 & 5 & 5 & 17623.048 \\
\hline 
\texttt{case6ww} & 6 & 3 & 11 & 4074.121 \\
\hline 
\texttt{pglib\_case14\_ieee} & 14 & 5 & 20 & 2665.655 \\
\hline 
\texttt{pglib\_case24\_ieee\_rts} & 24 & 33 & 38 & 52006.217 \\
\hline 
\texttt{pglib\_case30\_as} & 30 & 6 & 41 & 472.612 \\
\hline 
\texttt{pglib\_case30\_ieee} & 30 & 6 & 41 & 5485.280 \\
\hline 
\texttt{pglib\_case39\_epri} & 39 & 10 & 46 & 168303.267 \\
\hline 
\texttt{pglib\_case57\_ieee} & 57 & 7 & 80 & 33629.932 \\
\hline 
\texttt{pglib\_case73\_ieee\_rts} & 73 & 99 & 120 & 183995.780 \\
\hline 
\texttt{pglib\_case89\_pegase} & 89 & 12 & 210 & 109232.650 \\
\hline 
\texttt{pglib\_case118\_ieee} & 118 & 54 & 186 & 96881.525 \\
\hline 
\texttt{pglib\_case162\_ieee\_dtc} & 162 & 12 & 284 & 89234.092 \\
\hline 
\texttt{pglib\_case179\_goc} & 179 & 29 & 263 & 77442.263 \\
\hline 
\texttt{pglib\_case200\_activ} & 200 & 38 & 245 & 26231.327 \\
\hline 
\texttt{pglib\_case240\_pserc} & 240 & 143 & 448 & 49.12 \\
\hline 
\texttt{pglib\_case300\_ieee} & 300 & 69 & 411 & 521.98 \\
\hline
\end{tabular}
\caption{The attributes of the power test cases under consideration in the PG-lib library.}
\label{tab: Table1}
\end{table}
We employ a unique performance profile that measures CPU times to evaluate the efficiency of different solvers. The fundamental concept is as follows: for every test case $n$ and each solver $s$, we designate $t_{ns}$ as the time taken to solve the test case $n$ using solver $s$. We then establish the performance ratio as $\kappa_{ns}$. Assuming that $N_{tc}$ represents the total count of test cases analyzed, a comprehensive evaluation of solver $s$ performance for a specified $\tau$ is determined as follows. 
\begin{align}
P\left(\kappa_{ns} \leq \tau\right) & =\frac{\text{the count of test cases}}{N_{tc}} \\
\kappa_{ns} & =\frac{t_{ns}}{\min_s t_{ns}}
\end{align}
Fig. \ref{fig: fig11} presents the performance profile that details the CPU times for the IPOPT, \texttt{Baron 23.6.15}, and \texttt{Gurobi 10} methods applied to power system test cases. From our observations, it is clear that the proposed LCQP-ACOPF algorithm with bounded feasible SLP and \texttt{Gurobi 10} surpass the comparative \texttt{Baron 23.6.15} solver and the proposed Penalized-SDP ACOPF integrated with an unbounded SLP algorithm from \cite{47} in aspects of CPU time utilization and the number of test cases resolved. Specifically, Baron can only optimally solve 6 of the 20 test cases considered, with \texttt{case14-ieee} being the largest among them. The two methods, LCQP-ACOPF and Gurobi, prove to be more effective, as they manage to solve 12 test cases with the unbounded-SLP, 13 test cases using Gurobi, and 14 using LCQP-ACOPF within a 3-hour CPU time limit. Furthermore, the performance profile indicates that the LCQP-ACOPF method operates faster than Gurobi for these power test cases.

\begin{figure}[t]
\centerline{\includegraphics[width=0.8\textwidth]{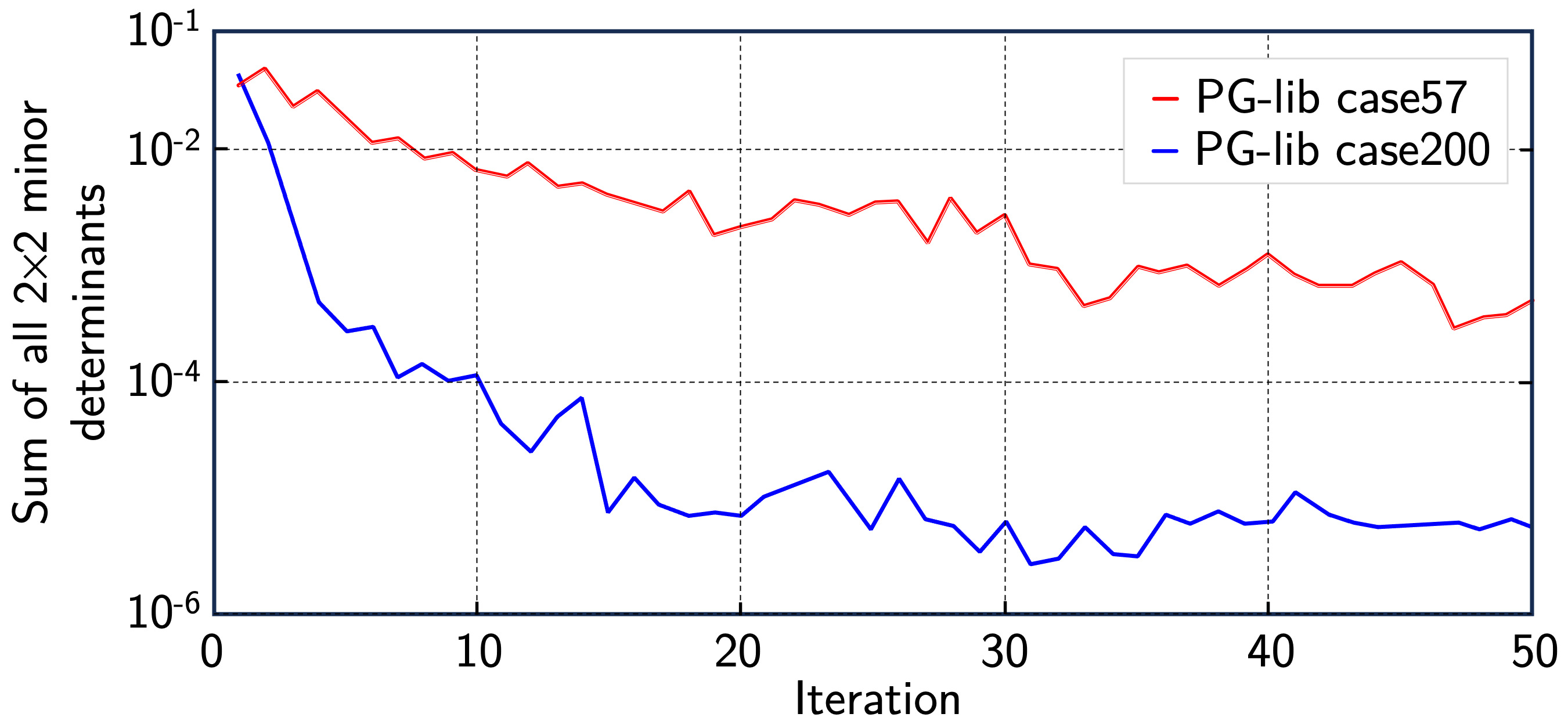}}
\caption{Total summation of all principal $2\times2$ minor determinants in the objective function over iterations of the 57-bus system and the 200-bus systems.}
\label{fig: fig10}
\end{figure}

For test cases involving 57 and 200 buses, we extended the iteration time to 50 and analyzed the total summation of all the principal $2\times2$ minor determinants in the objective function as a vanishing regularizer. The trend pertaining to the sum of all principal $2\times2$ minor determinants is illustrated in Fig. \ref{fig: fig10}. A notable exponential decline in the total sum of all principal minor determinants is visible within the first 25 iterations. This suggests that the $\operatorname{FBB}$ algorithm's solution could potentially attain feasibility for the lifted-ACOPF problem following a series of cut generations across the $(\Delta_1,\Delta_2)$ edges within the branch-and-bound algorithm.

{\small
\begin{table*}[t]
  \centering
\begin{tabular}{|l|c||c|c|c|c|c|}
\hline \multicolumn{2}{|c||}{ Test Case } & \multicolumn{5}{c|}{ LCQP-ACOPF } \\
\hline ID & Gap & Aux & UB & P-SDP & CPU & BB Nodes \\
\hline \texttt{pglib\_caseWB2} & 2.011 & 5 & 0 & 1 & 2 & 27  \\
\hline \texttt{pglib\_caseWB3} & 0.000 & 7 & 0 & 2 & 2 & 2  \\
\hline \texttt{pglib\_case3\_Imbd} & 1.670 & 7 & 0 & 2 & 1 & 2  \\
\hline \texttt{pglib\_caseWB5} & 14.481 & 13 & 0 & 2 & 360 & 8279  \\
\hline \texttt{pglib\_case5\_pjm} & 0.000 & 13 & 0 & 1 & 1 & 0  \\
\hline \texttt{pglib\_case6ww} & 0.000 & 14 & 5 & 1 & 6 & 0  \\
\hline \texttt{pglib\_case14\_ieee} & 0.100 & 32 & 0 & 2 & 1 & 0  \\
\hline \texttt{pglib\_case24\_ieee\_rts} & 0.000 & 52 & 0 & 2 & 1 & 1  \\
\hline \texttt{pglib\_case30\_as} & 0.050 & 62 & 8 & 0 & 8 & 0  \\
\hline \texttt{pglib\_case30\_ieee} & 0.000 & 62 & 0 & 1 & 1 & 0 \\
\hline \texttt{pglib\_case39\_epri} & 0.000 & 80 & 1 & 1 & 2 & 0  \\
\hline \texttt{pglib\_case57\_ieee} & 0.05 & 115 & 167 & 1 & - & 2659  \\
\hline \texttt{pglib\_case73\_ieee\_rts} & 0.000 & 146 & 67 & 2 & 70 & 0  \\
\hline \texttt{pglib\_case89\_pegase} & 0.000 & 180 & 31 & 3 & 34 & 0  \\
\hline \texttt{pglib\_case118\_ieee} & 2.102 & 238 & 933 & 2 & - & 973  \\
\hline \texttt{pglib\_case162\_ieee\_dtc} & 2.470 & 326 & 334 & 6 & - & 326  \\
\hline \texttt{pglib\_case179\_goc} & 0.024 & 360 & 1517 & 5 & - & 343  \\
\hline \texttt{pglib\_case200\_activ} & 0.000 & 402 & 423 & 5 & 428 & 0  \\
\hline \texttt{pglib\_case240\_pserc} & 2.470 & 482 & 3634 & 27 & - & 219  \\
\hline \texttt{pglib\_case300\_ieee} & 1.090 & 612 & 3278 & 42 & - & 34  \\
\hline
\end{tabular}
\caption{The initial gap, sizes of the power network, CPU times, and the number of nodes used in the branch-and-bound algorithm (FBB) for the LCQP-ACOPF.}
\label{tab: Table2}
\end{table*}
}

{\small
\begin{table*}[t]
  \centering
\begin{tabular}{|l|c||c|c|c|c|c|}
\hline \multicolumn{2}{|c||}{ Test Case } & \multicolumn{5}{c|}{ Unbounded-SLP in \cite{47} } \\
\hline ID & Gap & Aux & UB & P-SDP & CPU & BB Nodes \\
\hline \texttt{pglib\_caseWB2} & 2.011 & 10 & 0 & 1 & - & 44265 \\
\hline \texttt{pglib\_caseWB3} & 0.000 & 18 & 0 & 1 & 1 & 0 \\
\hline \texttt{pglib\_case3\_Imbd} & 1.670 & 23 & 0 & 1 & 1 & 0 \\
\hline \texttt{pglib\_caseWB5} & 14.481 & 38 & 0 & 1 & - & 292450 \\
\hline \texttt{pglib\_case5\_pjm} & 0.000 & 36 & 0 & 1 & 1 & 0 \\
\hline \texttt{pglib\_case6ww} & 0.000 & 60 & 0 & 2 & 2 & 0 \\
\hline \texttt{pglib\_case14\_ieee} & 0.100 & 110 & 0 & 1 & 4 & 0 \\
\hline \texttt{pglib\_case24\_ieee\_rts} & 0.000 & 186 & 56 & 1 & 55 & 0 \\
\hline \texttt{pglib\_case30\_as} & 0.050 & 226 & 0 & 1 & 1 & 0 \\
\hline \texttt{pglib\_case30\_ieee} & 0.000 & 230 & 9 & 1 & 8 & 0 \\
\hline \texttt{pglib\_case39\_epri} & 0.000 & 248 & 9 & 1 & 7 & 0 \\
\hline \texttt{pglib\_case57\_ieee} & 0.05 & 445 & 35 & 1 & - & 2278 \\
\hline \texttt{pglib\_case73\_ieee\_rts} & 0.000 & 580 & 191 & 3 & 142 & 0 \\
\hline \texttt{pglib\_case89\_pegase} & 0.000 & 1004 & 72 & 3 & 76 & 0 \\
\hline \texttt{pglib\_case118\_ieee} & 2.102 & 952 & 2866 & 5 & - & 1023 \\
\hline \texttt{pglib\_case162\_ieee\_dtc} & 2.470 & 1444 & 524 & 9 & - & 1045 \\
\hline \texttt{pglib\_case179\_goc} & 0.024 & 1246 & 3050 & 10 & - & 1132 \\
\hline \texttt{pglib\_case200\_activ} & 0.000 & 1284 & 2468 & 10 & 2176 & 0 \\
\hline \texttt{pglib\_case240\_pserc} & 2.470 & 1872 & 3330 & 29 & - & 819 \\
\hline \texttt{pglib\_case300\_ieee} & 1.090 & 2334 & 3743 & 52 & - & 3171 \\
\hline
\end{tabular}
\caption{The initial gap, sizes of the power network, CPU times, and the number of nodes used in the LCQP-ACOPF with unbounded-SLP methods, as presented in the ACOPF from \cite{47}.}
\label{tab: Table3}
\end{table*}
}

In Tables \ref{tab: Table2} and \ref{tab: Table3}, we conduct an exhaustive comparison between the LCQP-ACOPF and ACOPF methodologies, as well as the unbounded-SLP as delineated in \cite{47}. The ``Gap'' column is defined by the equation:

\begin{equation}
\text{Gap} = \left|\frac{\text{Opt} - \text{Cont}}{\text{Opt}}\right| \times 100
\end{equation}
This section describes key metrics and observations from an experimental analysis within the branch-and-bound framework for solving optimization problems, specifically focusing on ACOPF challenge using the LCQP model. The terminology used includes:

\begin{itemize}
    \item \textbf{Gap}: This term represents the initial difference encountered at the root node of the branch-and-bound process.
    \item \textbf{Cont}: Denotes the highest value achieved through rank relaxation.
    \item \textbf{Opt}: Defined according to Tables \ref{tab: Table2} and \ref{tab: Table3}, referring to the optimal solution.
    \item The notation ``-'' indicates that the algorithm failed to find a feasible solution.
    \item \textbf{Aux}: Shows the total number of auxiliary variables added in the relaxation process, equivalent to the total number of cuts applied in the spatial branch-and-bound method.
    \item \textbf{UB} and \textbf{P-SDP}: Document the time in seconds required to establish an initial upper bound and to solve the rank-one constraint relaxation on a standard CPU, respectively.
    \item \textbf{CPU}: Aggregates the total CPU time used, with ``-'' marking cases unsolved within a three-hour limit.
    \item \textbf{Nodes}: Reflects the number of nodes explored during the branch-and-bound process.
\end{itemize}

Key findings from this analysis include the following.

\begin{enumerate}
    \item The process of reformulating the problem is notably swift, often taking less than 428 seconds. This efficiency is attributed to the use of conventional semidefinite solvers to directly tackle rank relaxation.
    \item  In 9 out of 20 power test cases, there was no initial gap at the root node, highlighting the effectiveness of rank relaxation in addressing the OPF challenge. However, difficulties arise in precisely solving cases with a non-zero initial gap.
    \item The approach cited in \cite{47} had trouble solving smaller test cases with only 2 or 5 buses, due to its inability to improve the lower bound despite exploring many nodes.
    \item  Conversely, the LCQP-ACOPF method, while examining fewer nodes, was successful in progressively enhancing the lower bound for even the largest test cases. This success is largely due to the efficiency of the method in reducing the number of auxiliary variables and simplifying the relaxed equalities (\(s_i = t_i^2\)).
    \item Achieving a feasible solution can be challenging, particularly when using a generic interior point algorithm. The CPU time to find a feasible solution varies significantly, exemplified by the \texttt{case200\_activ} test case, where it ranged from 400 to 1380 seconds.
\end{enumerate}

These insights not only demonstrate the complexities and challenges involved in solving the ACOPF problem, but also underscore the potential of the LCQP-ACOPF method in achieving efficient and effective solutions under certain conditions.

\section{Conclusion}
The article explores the challenging LCQP model as it applies to the ACOPF problem, which involves complex variables with negative eigenvalues and is known for its NP-hard status, drawing from various fields. In our study, we have combined a series of advanced techniques: a feasible successive linear programming (FSLP) approach, robust convex relaxation methods, initial setup processes, and the branch-and-bound (B\&B) method to develop the FBB algorithm. This algorithm is specifically designed to find a globally optimal solution for the fundamental LCQP-based ACOPF model. Furthermore, we have developed a new joint algorithm that operates independently of the B\&B framework. This innovative algorithm aims to discover a global optimum for the LCQP-ACOPF problem by integrating FSLP, convex relaxation, and a unique B\&B approach. With certain basic assumptions, we have verified the global effectiveness of both the FBB and FSLP algorithms and estimated their computational complexity. Early experimental results confirm the success of FSLP and FBB in consistently finding the global optimum for all the power system scenarios described in our simulations. It is crucial to note that the complexity of the FBB algorithm increases exponentially with the number of negative eigenvalues in the Hessian matrix of the objective function, suggesting potential inefficiencies for power test cases with many negative eigenvalues. Our research shows that a specific implementation of our approach allows for the use of carefully penalized semidefinite relaxation, leading to much tighter convex relaxations compared to traditional methods. This significantly improves our handling of medium-sized real-world power test cases related to the LCQP-ACOPF problem. Future research could fruitful explore advanced global algorithms for the broader ACOPF problem by integrating FSLP with additional branching and bound or local search strategies. Another promising area of investigation involves creating similar global solutions for generic quadratic programming problems with non-convex constraints.

\appendix




\end{document}